\documentclass[twoside,leqno,11pt]{article}
\usepackage{revstat}
\usepackage{amsmath}
\usepackage{graphicx,psfrag,epsf}
\usepackage{enumerate}
\usepackage{amsfonts}
\usepackage{amsthm}
\usepackage{mathrsfs}
\usepackage{multirow}
\usepackage[titletoc,title]{appendix}
\usepackage[font={footnotesize}]{caption}
\usepackage[font={footnotesize}]{subcaption}

\newtheorem{thm}{Theorem}[section]

\newtheorem{cor}[thm]{Corollary}

\newtheorem{rem}{Remark}[section]

\newcommand\ba{\begin{array}}
	\newcommand\ds{\displaystyle}
	\newcommand\ms{\medskip}
	\newcommand\ea{\end{array}}

\newcommand*\samethanks[1][\value{footnote}]{\footnotemark[#1]}

\begin{document}
	
	\bibliographystyle{Chicago}

\title{INFERENCE FOR MULTIVARIATE REGRESSION MODEL BASED ON SYNTHETIC DATA\\
GENERATED UNDER FIXED-POSTERIOR\\
PREDICTIVE SAMPLING: COMPARISON WITH\\
PLUG-IN SAMPLING}
	\renewcommand{\titleheading}{Inference for MLR model under FPPS: Comparison with Plug-in Sampling}
	
	\author{\authoraddress{Ricardo Moura}
		{CMA,
			Faculty of Sciences and Technology,
			Nova University of Lisbon\\
			Portugal
			\ (rp.moura@campus.fct.unl.pt)}
		\\
		\authoraddress{Martin Klein\thanks{
				\textbf{Disclaimer:} This article is released to inform interested parties of ongoing research and to encourage discussion. The views 
				expressed are those of the authors and not necessarily those of the U.S. Census Bureau.}}
		{Center for Statistical Research and Methodology,
			U.S. Census Bureau,\\
			U.S.A.
			\ (martin.klein@census.gov)}
		\\
		\authoraddress{Carlos A. Coelho}
		{CMA and Mathematics Department, 
			Faculty of Sciences and Technology,\\ Nova University of Lisbon\\
			Portugal
			\ (cmac@fct.unl.pt)}
		\\
		\authoraddress{Bimal Sinha\samethanks}
		{Department of Mathematics and Statistics,\\
				University of Maryland, Baltimore County\\ 
			and Center for Disclosure Avoidance Research,
			U.S. Census Bureau\\
			U.S.A.
			\ (sinha@umbc.edu)}
	}
	\renewcommand{\authorheading}
	{\ Ricardo Moura, \ Martin Klein, \ Carlos A. Coelho and \ Bimal Sinha}  

	\date{\vspace{-5ex}}
	
	\maketitle

\vspace{-.65cm}	
	
\begin{abstract} 
		The authors derive likelihood-based exact inference methods for the multivariate regression
model, for singly imputed synthetic data
generated via Posterior Predictive Sampling (PPS) and for multiply imputed synthetic data generated via a newly proposed sampling method, which the authors call Fixed-Posterior Predictive Sampling (FPPS). In the single imputation case, our proposed FPPS method concurs with the
usual Posterior Predictive Sampling (PPS) method, thus filling the gap in the existing literature where inferential methods are only available for multiple imputation. Simulation studies compare the results obtained with those for the exact test procedures under the Plug-in Sampling method, obtained by the same authors. Measures of privacy are discussed and compared with the measures derived for the Plug-in Sampling method. An application using U.S.\ 2000 Current Population Survey data is discussed.
\end{abstract}
		
\vspace{-.5cm}
		


\vspace{-.5cm}

\begin{ams}
		62H10, 62H15, 62H12, 62J05, 62F10, 62E15, 62E10, 62E17, 62D99.
\end{ams}

\mainpaper
	
\section{INTRODUCTION}
	
	When releasing microdata to the public, methods of statistical disclosure control (SDC) are used to protect confidential data, that is ``data which allow statistical units to be identified, either directly or indirectly, thereby disclosing individual information'' \cite{REG},
while enabling valid statistical inference to be drawn on the relevant population. SDC methods include data swapping, additive and multiplicative noise, top and bottom coding, and also the creation of synthetic data. In this paper, the authors provide inferential tools for the statistical analysis of a singly imputed synthetic dataset when the real dataset cannot be released. The multiple imputation case is also addressed, using a new adapted method of generating synthetic data, which the authors call Fixed-Posterior Predictive Sampling (FPPS).
	
	The use of synthetic data for SDC started with Little \cite{little93} and Rubin \cite{rubin93} using multiple imputation \cite{rubin87}. Reiter \cite{reiter03} was the first to present methods for drawing inference based on partially synthetic data. Moura et al. \cite{moura16} complemented this work with the development of likelihood-based exact inference methods for both single and multiple imputation, that is, inferential procedures developed based on exact distributions, and not on asymptotic results, in the case where synthetic datasets were generated via Plug-in Sampling. The procedures of Reiter \cite{reiter03} are general in that they can be applied to a variety of estimators and statistical models, but these procedures are only applicable in the multiple imputation case, and are based on large sample approximations.
	
	There are two major objectives in the present research. First, to make available likelihood-based exact inference for singly imputed synthetic data via Posterior Predictive Sampling (PPS) where the usual available procedures are not applicable, therefore extending the work of Klein and Sinha \cite{Klein2015}, under the multivariate linear regression (MLR) model. Second, to propose a different approach for release of multiple synthetic datasets, FPPS, which can use a similar way of gathering information from the synthetic datasets to that used in \cite{moura16}, when these synthetic datasets are generated via the Plug-in Sampling method. This second objective arises from the fact that when using the classical PPS it is too hard to construct an exact joint probability density function (pdf) for the estimators, under the MLR model, since one would face the problem of deriving the distribution of a sum of variables that follow Wishart distributions with different parameter matrices. It is with this problem in mind, that we propose an adapted method that we will call the FPPS method. We show that this method offers a higher level of confidentiality than the Plug-in Sampling method, and it still allows one to draw inference for the unknown parameters using a joint pdf of the proposed estimators.
	
	A brief description of the PPS and FPPS methods follows. Suppose that 
	$\mathbf{Y}=(\mathbf{y}_1,...,\mathbf{y}_n)$ are 
	the original data which are jointly distributed according to the pdf 
	$f_{\boldsymbol\theta}(\mathbf{Y})$, where $\boldsymbol\theta$ is the unknown (scalar, vector or matrix) parameter. A prior $\pi(\mathbf{\boldsymbol{\theta}})$ for $\mathbf{\boldsymbol{\theta}}$ is assumed and then the posterior distribution of $\mathbf{\boldsymbol{\theta}}$ is obtained as $\pi(\boldsymbol{\theta}|Y)\propto \pi(\mathbf{\boldsymbol{\theta}})f_{\mathbf{\boldsymbol{\theta}}(x)}$, and used to draw a replication $\mathbf{\boldsymbol{\theta}}_f^\bullet$ of $\mathbf{\boldsymbol{\theta}}$, when applying the FPPS, or draw $M \geq 1$ independent replications $\mathbf{\boldsymbol{\theta}}_1^\bullet,...,\mathbf{\boldsymbol{\theta}}_M^\bullet$ of $\boldsymbol{\theta}$, when applying the PPS.
	In the case of FPPS, we generate $M$ replicates of $\mathbf{Y}$, namely, $\mathbf{W}_j=(\mathbf{w}_{j1},...,\mathbf{w}_{jn})$, $j=1,...,M$ drawn all independently from the same $f_{\boldsymbol{\theta}_f^\bullet}$, where $f_{\boldsymbol{\theta}_f^\bullet}$ is the joint pdf of the original $\mathbf{Y}$ with $\boldsymbol{\theta}_f^\bullet$ replacing the unknown $\boldsymbol{\theta}$. In the case of the usual PPS method for each $j$-th generated synthetic dataset we would use the corresponding $j$-th posterior draw $\boldsymbol{\theta}_j^\bullet$ and corresponding $j$-th joint pdf's $f_{\boldsymbol{\theta}_j^\bullet}$, for $j=1,...,M$. In either case, these synthetic datasets $\mathbf{W}_1, \hdots, \mathbf{W}_M$ will be the datasets available to the general public. One may observe that, for $M=1$, the Posterior Predictive Sampling and Fixed-Posterior Predictive Sampling methods concur.

	Regarding the MLR model, in our context, we consider the \textit{sensitive} response variables $y_j$ $(j=1,...,m)$ forming the vector of response variables $\mathbf{y}={(y_1,...,y_m)'}$, and a set of p non-\textit{sensitive} explanatory variables $\mathbf{x}=(x_1,...,x_p)'$. It is assumed that $\mathbf{y}|\mathbf{x}\sim N_m(\mathbf{B'}\mathbf{x},\mathbf{\Sigma})$, with $\mathbf{B}$ and $\mathbf{\Sigma}$ unknown, and the original data consist of $\mathcal{Y}=\linebreak\lbrace(y_{1i},...,y_{mi},x_{1i},...,x_{pi}),i=1,...,n\rbrace$, where $n$ will be the sample size. Let us consider $\mathbf{Y}=(\mathbf{y}_1,...,\mathbf{y}_n)$ with $\mathbf{y}_i=(y_{1i},...,y_{mi})'$ and $\mathbf{X}=(\mathbf{x}_1,...,\mathbf{x}_n)$ with $\mathbf{x}_i=(x_{1i},...,x_{pi})'$. We assume $rank(\mathbf{X}:p \times n) = p < n$ and $n\geq m+p$. Therefore the following regression model is considered
	\begin{equation}\label{eq:model}
	\mathbf{Y}_{m\times n}=\mathbf{B}'_{m\times p} \mathbf{X}_{p\times n}  + \mathbb{E}_{m\times n},
	\end{equation}
	\noindent where $\mathbb{E}_{m\times n}$ is distributed as $N_{mn}(\mathbf{0},\mathbf{I}_n\otimes \mathbf{\Sigma})$. Based on the 
	original data, 
	\begin{equation}\label{eq:orB}
	\mathbf{\hat{B}}=(\mathbf{XX}')^{-1}\mathbf{XY'}
	\end{equation}
	is the Maximum Likelihood Estimator (MLE) and the Uniformly Minimum-Variance Unbiased Estimator (UMVUE) of $\mathbf{B}$, distributed as 
	$N_{pm}(\mathbf{B},\mathbf{\Sigma}\otimes(\mathbf{XX'})^{-1})$, independent of 
	$\mathbf{\hat{\Sigma}}=\frac{1}{n}(\mathbf{Y}-\mathbf{\hat{B}'}\mathbf{X})(\mathbf{Y}-\mathbf{\hat{B}'}\mathbf{X})'$ which is the MLE 
	of $\mathbf{\Sigma}$, with $n \mathbf{\hat{\Sigma}}\sim W_m(\mathbf{\Sigma},n-p)$. Therefore
	\begin{equation}\label{eq:orS}
	\mathbf{S}=\frac{n\mathbf{\hat{\Sigma}}}{n-p}
	\end{equation}
	will be the UMVUE of $\mathbf{\Sigma}$.

	The organization of the paper is as follows. In Section \ref{sec:post}, based on singly and multiply imputed synthetic datasets generated via Fixed-Posterior Predictive Sampling, two procedures are proposed to draw inference for the matrix of regression coefficients. Under the single imputation case, we recall that the FPPS and the PPS methods coincide. The test statistics proposed will be pivot statistics, different from the classical test statistics for $\mathbf{B}$ under the MLR model (see \cite[Secs 8.3 and 8.6]{anderson84}) since it is shown that these classical test statistics are not pivotal in the present context. Section \ref{sec:sim} presents some simulations in order to check the accuracy of theoretically derived results. Also in this section, the authors use a measure for the \textit{radius} (distance between the center and the edge) of the confidence sets for the regression coefficients adapted from \cite{moura16}, computed for the original data and also for the synthetic data generated via FPPS. These \textit{radius} measures are compared with the ones obtained when synthetic datasets are generated via Plug-in Sampling. Section \ref{sec:app} presents data analyses under the proposed methods in the context of public use data from the U.S. Current Population Survey comparing with the same data analysis given by \cite{moura16} under the Plug-in Sampling method. In Section \ref{sec:pri}, we compare the level of privacy protection obtained via our FPPS method and via Plug-in Sampling method. Some concluding remarks are added in Section \ref{sec:con}. Proofs of the theorems, and other technical derivations are presented in Appendices \ref{App:A} and \ref{App:last}.

\section{ANALYSIS FOR SINGLE AND MULTIPLE IMPUTATION}
\label{sec:post}
	
	In this section, we present two new exact likelihood-based procedures for the analysis of synthetic data generated using Fixed-Posterior Predictive Sampling method, under the MLR model in (\ref{eq:model}). For the single imputation case, the two new procedures developed also offer the possibility of drawing inference for a single synthetic dataset generated via Posterior Predictive Sampling.

\subsection{A FIRST NEW PROCEDURE}
\label{ssec:mul1}
	
	In this subsection, the synthetic data will consist of $M$ synthetic versions of $\mathbf{Y}$ generated based on the FPPS method. 
	
	Consider the joint prior distribution $\pi(\mathbf{B},\mathbf{\Sigma})\propto |\mathbf{\Sigma}|^{-\alpha/2},$ leading to the posterior distributions for $\mathbf{\Sigma}$ and $\mathbf{B}$
	\begin{equation}\label{eq:postSigma}
	\mathbf{\Sigma}|_{\mathcal{Y},\mathbf{S}}\sim W^{-1}_m ((n-p)\mathbf{S},n+\alpha-p)
	\end{equation} and
	\begin{equation}\label{eq:postB}
	\mathbf{B}|_{\mathcal{Y},\mathbf{\Sigma}}\sim N_{pm}(\mathbf{\hat{B},\Sigma\otimes(XX')^{-1}}),
	\end{equation}
	where we assume that $n+\alpha>p+m+1$ (see proof in Appendix \ref{Aapp:last0}). Consequently, we draw $\mathbf{\tilde{\Sigma}}$ from (\ref{eq:postSigma}) and $\mathbf{\tilde{B}}$ from (\ref{eq:postB}), upon replacing $\mathbf{\Sigma}$ by $\mathbf{\tilde{\Sigma}}$ in this latter expression. We then generate the $M$ synthetic datasets, denoted as  $\mathbf{W}_j=(\mathbf{w}_{j1},...,\mathbf{w}_{jn})$, for $j=1,...,M$, where $\mathbf{w}_{ji}=(w_{1ji},...,w_{mji})'$, are independently distributed as
	\begin{equation}\label{eq:synt}
	\mathbf{w}_{ji}|_{\mathbf{\tilde{B}},\mathbf{\tilde{\Sigma}}}\sim N_m(\mathbf{\tilde{B}'x_i},\mathbf{\tilde{\Sigma}}),~~ i=1,...,n, j=1,...,M.
	\end{equation}
	For $i=1,...,n$ and $j=1,...,M$, let $\mathbf{B}_j^{\bullet}=(\mathbf{XX'})^{-1}\mathbf{XW}'_j $ and $\mathbf{S}_j^{\bullet}=\frac{1}{n-p}(\mathbf{W}_j-\mathbf{B}_j^{\bullet'}\mathbf{X})(\mathbf{W}_j-\mathbf{B}_j^{\bullet'}\mathbf{X})'$
	be the estimators of $\mathbf{B}$ and $\mathbf{\Sigma}$, based on the synthetic data ${(w_{1ji},...,w_{mji},x_{1i},...,x_{pi})}$, which by Lemma 1.1 in \cite{moura16} are jointly sufficient. Conditional on $(\mathbf{\tilde{B}},\mathbf{\tilde{\Sigma}})$, for every $j=1,...,M$, $\mathbf{B}_j^{\bullet}$ is independent of $\mathbf{S}_j^{\bullet}$ and $\left\{(\mathbf{B}_1^{\bullet},\mathbf{S}_1^{\bullet}),...,(\mathbf{B}_M^{\bullet},\mathbf{S}_M^{\bullet})\right\}$ are jointly sufficient estimators for $\mathbf{B}$ and $\mathbf{\Sigma}$. Define then 
	\begin{equation}\label{eq:parameters1st}
	\overline{\mathbf{B}}^\bullet_M=\frac{1}{M}\sum_{j=1}^M \mathbf{B}_j^{\bullet}~~~{\rm and}~~~\overline{\mathbf{S}}^\bullet_M=\frac{1}{M}\sum_{j=1}^M \mathbf{S}_j^{\bullet},
	\end{equation}
	which are also mutually independent, given $\tilde{\mathbf{B}}$ and $\tilde{\mathbf{\Sigma}}$.
	For $p\geq m$ and
	$n+\alpha>p+2m+2$, we derive the following main results.
	\vspace{-10pt}
	\begin{enumerate}
		
		\item The MLE of $\mathbf{B}$ is $\overline{\mathbf{B}}^\bullet_M$, which is unbiased for $\mathbf{B}$, with $Var(\overline{\mathbf{B}}^\bullet_M)\!$\linebreak$=\!N_{M,n,m,p,\alpha}\mathbf{\Sigma}\!\otimes\!(\mathbf{\!XX'})^{-1}$, where 
		$N_{M,n,m,p,\alpha}=\frac{2M(n+\frac{\alpha}{2}-p-m-1)+n-p}{M(n+\alpha-p-2m-2)}$
		(see Theorem \ref{thm:pdf} and Appendix \ref{Aapp:last3}).
		
		\item An unbiased estimator (UE) of $\mathbf{\Sigma}$ will be $\mathbf{\hat{S}}_M=\frac{n+\alpha-p-2m-2}{n-p}\mathbf{\overline{S}}^{\bullet}_M$ (see Theorem \ref{thm:pdf} and Appendix \ref{Aapp:last3}); for $\alpha=2m+2$, $\mathbf{\overline{S}}^{\bullet}_M$ will also be an UE for $\mathbf{\Sigma}$,
		
		\item In Theorem \ref{thm:dist} (see below), we prove that 
		\begin{equation}\label{eq:T1st}
		T^\bullet_M=\frac{|(\overline{\mathbf{B}}^{\bullet}_M-\mathbf{B})'(\mathbf{XX}')(\mathbf{\overline{B}}^{\bullet}_M-\mathbf{B)}|}{|M(n-p)\mathbf{\overline{S}}^{\bullet}_M|}\,,
		\end{equation}
		a statistic somewhat related with the Hotelling $T^2$, this one built to make inference on a matrix parameter, is a pivotal quantity, and that for $\mathbf{A}_1\sim W_m(\mathbf{I}_m,n+\alpha-p-m-1)$, $\mathbf{A}_2\sim W_m(\mathbf{I}_m,n-p)$ and $F_i\sim F_{p-i+1,M(n-p)-i+1}$\linebreak $(i=1,...,m)$, all independent random variables,  
		$$T^\bullet_M|_{\mathbf{\Omega}}\stackrel{st}{\sim}\left\{\prod_{i=1}^{m}\frac{p-i+1}{M(n-p)-i+1}F_{i}\right\}\left|\frac{M+1}{M}\mathbf{I}_{m}+\mathbf{\Omega}\right|,$$
		where $\mathbf{\Omega}$ has the same distribution as $\mathbf{A}_1^{\frac{1}{2}}\mathbf{A}_2^{-1}\mathbf{A}_1^{\frac{1}{2}}$ and where $\stackrel{st}{\sim}$ means `stochastic equivalent to'.

		\item If one wants to test a linear combination of the parameters in $\mathbf{B}$, namely, $\mathbf{C = A B}$ where $\mathbf{A}$ is a $k\times p$ matrix with $rank(\mathbf{A})=k\leq p$ and $k\geq m$, one defines 
		$$T^\bullet_{M,\mathbf{C}}=\linebreak\frac{|(\mathbf{A \overline{B}}^{\bullet}_M-\mathbf{C})'(\mathbf{A(XX')^{-1}A'})^{-1}(\mathbf{A\overline{B}}^{\bullet}_M-\mathbf{C})|}{|M(n-p)\mathbf{\overline{S}}^{\bullet}_M|}$$and proceeds by noting that
		\begin{equation}\label{eq:TC1}
		T^\bullet_{M,\mathbf{C}}|_\mathbf{W}\stackrel{st}{\sim}\left\{ \prod_{i=1}^{m}\frac{k-i+1}{M(n-p)-i+1}F_{k,i}\right\}\left|\frac{M+1}{M}\mathbf{I}_{m}+\mathbf{\Omega}\right|,
		\end{equation}
		with $F_{k,i}\sim F_{k-i+1,M(n-p)-i+1}$ being independent random variables and $\mathbf{\Omega}$ defined as in the previous item.
		
		(i)\textit{Test for the significance of $\mathbf{C}$:} in order to test $H_0:\mathbf{C}=\mathbf{C}_0$ versus 
		$H_1:\mathbf{C}\neq\mathbf{C}_0$, we reject $H_0$ whenever $T^\bullet_{M,\mathbf{C}_0}$ exceeds $\delta_{M,k,m,p,n;\gamma}$ where 
		$\delta_{M,k,m,p,n;\gamma}$ satisfies $(1-\gamma)= Pr(T^\bullet_{M,\mathbf{C}_0}\leq\delta_{M,k,m,p,n;\gamma})$ when $H_0$ is true. To perform a test for $\mathbf{B=B}_0$ one has to take $\mathbf{A=I}_p$.
		
		(ii)\textit{Confidence set for $\mathbf{C}$:} a $(1-\gamma)$ level confidence set for $\mathbf{C}$ is given by
		\begin{equation}\label{eq:elipsoidAB}
		\Delta_M(\mathbf{C})=\lbrace\mathbf{C}:T^\bullet_{M,\mathbf{C}}\leq\delta_{M,k,m,n,p;\gamma}\rbrace,
		\end{equation}
		where the value of $\delta_{M,k,m,n,p;\gamma}$ can be obtained by simulating the distribution in (\ref{eq:TC1}).
	\end{enumerate}
	Results in 1-4 are derived based on Theorems \ref{thm:pdf} and \ref{thm:dist} below.
	
	\begin{thm}\label{thm:pdf}
		
		The joint pdf of $\mathbf{\overline{B}}^{\bullet}_M,\mathbf{\overline{S}}^{\bullet}_M$ and $\mathbf{\tilde{\Sigma}^{-1}}$, for $\mathbf{\overline{B}}^{\bullet}_M$ and $\mathbf{\overline{S}}^{\bullet}_M$ defined in (\ref{eq:parameters1st}), is proportional to
		\[\begin{array}{l}
		e^{-\frac{1}{2}tr\lbrace\mathbf{(\frac{M+1}{M}\tilde{\Sigma}+\Sigma)}^{-1}(\mathbf{\overline{B}}^{\bullet}_M-\mathbf{B})'\mathbf{XX'}(\mathbf{\overline{B}}^{\bullet}_M-\mathbf{B})+M(n-p)\mathbf{\tilde{\Sigma}^{-1}}\mathbf{\overline{S}}^{\bullet}_M\rbrace}\\ 
		~~~~~\times \frac{\vert\mathbf{\overline{S}}^{\bullet}_M\vert^{\frac{M(n-p)-m-1}{2}}}{\vert\mathbf{\tilde{\Sigma}}\vert^{\frac{M(n-p)+n+\alpha}{2}-m-1}}\vert\mathbf{\Sigma}\vert^{-\frac{n}{2}} |\frac{M}{M+1}\mathbf{\tilde{\Sigma}^{-1}}+\mathbf{\Sigma^{-1}}|^{-p/2} |\mathbf{\tilde{\Sigma}^{-1}+\Sigma^{-1}}|^{-\frac{2n+\alpha-2p-m-1}{2}},
		\end{array}\]
		so that $\overline{\mathbf{B}}^\bullet_M$ and $\overline{\mathbf{S}}^\bullet_M$, given $\tilde{\mathbf{\Sigma}}$, are independent, with $$\mathbf{\overline{B}}^{\bullet}_M|_{\tilde{\mathbf{\Sigma}}}\sim N_{pm}\left(\mathbf{B},\left(\frac{M+1}{M}\mathbf{\tilde{\Sigma}+\Sigma}\right)\otimes(\mathbf{XX}')^{-1}\right)$$
		and 
		$$\mathbf{\overline{S}}^{\bullet}_M|_{\tilde{\mathbf{\Sigma}}}\sim W_m\left(\frac{1}{M(n-p)}\mathbf{\tilde{\Sigma}},M(n-p)\right).$$

	\end{thm}
	
	\begin{proof}
		See Appendix \ref{App:A}.
	\end{proof}
	
	\begin{thm}\label{thm:dist}
		The distribution of the statistic $T^\bullet_M$ defined in (\ref{eq:T1st}) can be obtained from the decomposition
		$$T^\bullet_M|_{\mathbf{\Omega}}\stackrel{st}{\sim}\left\{\prod_{i=1}^{m}\frac{p-i+1}{M(n-p)-i+1}F_{i}\right\}\left\vert\frac{M+1}{M}\mathbf{I}_{m}+\mathbf{\Omega}\right\vert$$
		where $F_i\sim F_{p-i+1,M(n-p)-i+1}$ are independent random variables, themselves independent of $\mathbf{\Omega}$, which has the same distribution as $\mathbf{A}_1^{\frac{1}{2}}\mathbf{A}_2^{-1}\mathbf{A}_1^{\frac{1}{2}}$ with $\mathbf{A}_1\sim W_m(\mathbf{I}_m,n+\alpha-p-m-1)$ and $\mathbf{A}_2\sim W_m(\mathbf{I}_m,n-p)$, two independent random variables.
	\end{thm}
	
	\begin{proof}
		See Appendix \ref{App:A}.
	\end{proof}
	
	\begin{rem}
		\indent When $m=1$ and $M=1$, the statistic in (\ref{eq:T1st}) reduces to the statistic $T^2$ used in \cite{Klein2015} whose pdf is obtained by noting that
		$$T^2|_{\Omega=\omega}\sim \frac{p}{n-p}(2+\omega)F_{p,n-p} ~~~{\rm where}~~~ f_{\Omega}(\omega)\propto\frac{\omega^{\frac{n+\alpha-p-4}{2}}}{(1+\omega)^{\frac{2n+\alpha-2p-2}{2}}}.$$
	\end{rem}

		\begin{rem}
			We remark that the statistic $T^\bullet_M$ in (\ref{eq:T1st}) degenerates towards zero when $n\rightarrow\infty$ or $M\rightarrow\infty$, but 
			$$
			(M(n-p))^m\,T^\bullet_M|_{\mathbf{\Omega}}\xrightarrow[n\rightarrow\infty]{d}\left\{\prod_{i=1}^{m}\chi^2_{p-i+1}\right\}\left|\frac{M+1}{M}\mathbf{I}_{m}+\mathbf{\Omega}\right|
			$$
			and
			$$
			(M(n-p))^m\,T^\bullet_M|_{\mathbf{\Omega}}\xrightarrow[M\rightarrow\infty]{d}\left\{\prod_{i=1}^{m}\chi^2_{p-i+1}\right\}\left|\mathbf{I}_{m}+\mathbf{\Omega}\right|,
			$$
			where $\xrightarrow{~~d~~}$ represents convergence in distribution.
			Consequently, if instead of using $T^\bullet_M$ one uses $T^\bullet_{M2}=(M(n-p))^m\,T^\bullet_M=\frac{|(\overline{\mathbf{B}}^{\bullet}_M-\mathbf{B})'(XX')(\mathbf{\overline{B}}^{\bullet}_M-\mathbf{B)}|}{|\mathbf{\overline{S}}^{\bullet}_M|}$
			one would have
			$$
			T^\bullet_{M2}|_{\mathbf{\Omega}}\xrightarrow[n\rightarrow\infty]{d}\left\{\prod_{i=1}^{m}\chi^2_{p-i+1}\right\}\left|\frac{M+1}{M}\mathbf{I}_{m}+\mathbf{\Omega}\right|
			$$
			and
			$$
			T^\bullet_{M2}|_{\mathbf{\Omega}}\xrightarrow[M\rightarrow\infty]{d}\left\{\prod_{i=1}^{m}\chi^2_{p-i+1}\right\}\left|\mathbf{I}_{m}+\mathbf{\Omega}\right|,
			$$
			which corresponds to the use of a simple scale change.
		\end{rem}
	
	In Table 1, we list the simulated $0.05$ cut-off points for $T^\bullet_M$, for $M=1$ for some values of $p$, $m$ and $n$.
	
	\begin{table}[h!]
		\footnotesize
		\caption{\footnotesize{Cut-off points of the 95\% confidence set for the regression coefficient $\mathbf{B}$}}
		\vspace{-15pt}
		\begin{center}
			\begin{tabular}{c}
				\begin{tabular}{  | r | l | l | l | l | }
					\hline
					\multirow{3}{*}{$n$} & \multicolumn{4}{| c |}{$p=3$}\\[-2pt]
					& \multicolumn{2}{c|}{$m=1$} & \multicolumn{2}{c|}{$m=3$}\\[-2pt]
					&	\multicolumn{1}{c|}{$\alpha=2$} &	\multicolumn{1}{c|}{$\alpha=4$} &	\multicolumn{1}{c|}{$\alpha=4$} & 	\multicolumn{1}{c|}{$\alpha=6$}    \\ \cline{1-5}
					10 & 6.568     & 7.433     & 20.11     & 29.08\\
					50                    & 5.502E-01 & 5.581E-01 & 9.277E-03 & 9.691E-03\\ 
					100                   & 2.518E-01 & 2.542E-01 & 9.212E-04 & 9.443E-04\\ 
					200                   & 1.207E-01 & 1.208E-01 & 1.049E-04 & 1.064E-04\\ \hline
				\end{tabular}
				\\
				\begin{tabular}{ | r | l | l | l | l |}
					\hline
					\multirow{3}{*}{$n$} &\multicolumn{4}{| c |}{$p=4$}\\[-2pt]
					&\multicolumn{2}{|c|}{$m=1$} & \multicolumn{2}{c|}{$m=3$}\\[-2pt]
					&\multicolumn{1}{|c|}{$\alpha=2$} &	\multicolumn{1}{c|}{$\alpha=4$} &	\multicolumn{1}{c|}{$\alpha=4$} & 	\multicolumn{1}{c|}{$\alpha=6$}    \\ \cline{1-5}
					10 & 11.08      & 12.69     & 239.2      & 372.7 \\	 
					50 & 6.884E-01  & 6.984E-01 & 3.550E-02  & 3.697E-02\\ 
					100 & 3.108E-01  & 3.128E-01 & 3.487E-03  & 3.564E-03\\ 
					200 & 1.487E-01  & 1.490E-01 & 3.674E-04  & 3.723E-04 \\ \hline
				\end{tabular}
			\end{tabular}
		\end{center}
		\vspace{-8pt}
	\end{table}
	
	Similar to what was done in \cite{moura16}, one could suggest the following adaptations of the classical test criterion for the multivariate regression model (see \cite[Secs 8.3 and 8.6]{anderson84} for the classical criteria):
	\begin{itemize} 
		
		\item[(a)] $T^\bullet_{1,M}=|\mathbf{\overline{S}}_M^{\bullet}||\mathbf{\overline{S}}_M^{\bullet}+(\mathbf{\overline{B}}_M^{\bullet}-\mathbf{B})'(XX')(\mathbf{\overline{B}}_M^{\bullet}-\mathbf{B})|^{-1}$ (Wilks' Lambda 
		Criterion),
		\item[(b)] $T^\bullet_{2,M}=tr\left[(\mathbf{\overline{B}}^{\bullet}_M-\mathbf{B})'(\mathbf{XX'})(\mathbf{\overline{B}}^{\bullet}-\mathbf{B})(\mathbf{\overline{S}}^{\bullet}_M)^{-1}\right]$ (Pillai's Trace Criterion),
		\item[(c)] 
		$T^\bullet_{3,M}=tr\left[(\mathbf{\overline{B}}^{\bullet}_M\!-\!\mathbf{B})'(\mathbf{XX'})(\mathbf{\overline{B}}^{\bullet}_M\!-\!\mathbf{B})[(\mathbf{\overline{B}}^{\bullet}_M\!-\!\mathbf{B})'(\mathbf{XX'})(\mathbf{\overline{B}}^{\bullet}_M\!-\!\mathbf{B})+\mathbf{\overline{S}}^{\bullet}_M]^{-1}\!\right]$ (Hotelling-Lawley Trace Criterion),
		\item[(d)] $T^\bullet_{4,M}=\lambda_1$ where $\lambda_1$ denotes the largest eigenvalue of $(\mathbf{\overline{B}}_M^{\bullet}-\mathbf{B})'(\mathbf{XX'})(\mathbf{\overline{B}}_M^{\bullet}-\mathbf{B})(\mathbf{\overline{S}}^{\bullet}_M)^{-1}$ (Roy's 
		Largest Root Criterion).
		
	\end{itemize}
	
		However, these statistics are non-pivotal, since their distributions are function of $\mathbf{\Sigma}$ (see Appendix \ref{Aapp:last3}).

\subsection{A SECOND NEW PROCEDURE}
\label{ssec:mul2}
	
	We propose yet another likelihood-based approach for exact inference about $\mathbf{B}$ where one may gather more information from the released synthetic data, following a somewhat similar procedure to the one used in \cite{moura16}. Let us start by recalling that $\mathbf{W_j}\;(j=1,...,M)$ are $m\times n$ matrices formed by the vectors $(\mathbf{w}_{j1},...,\mathbf{w}_{jn})$ as columns, generated from $\mathbf{w}_{ji}|_{\mathbf{\tilde{B}},\mathbf{\tilde{\Sigma}}}\sim N_m(\mathbf{\tilde{B}'x_i},\mathbf{\tilde{\Sigma}})$ $(i=1,...,n)$.
	Note that, conditionally on $\mathbf{\tilde{B}}$ and $\mathbf{\tilde{\Sigma}}$, $(\mathbf{w}_{1i},...,\mathbf{w}_{Mi})$ is a random sample from $N_m(\mathbf{\tilde{B}'x}_i,\mathbf{\tilde{\Sigma}})$, for $i=1,...,n$. Consider $\mathbf{\overline{w}}_i=\frac{1}{M}\sum_{j=1}^M\mathbf{w}_{ji}$ and $\mathbf{S}_{\mathbf{w}i}=\sum_{j=1}^M (\mathbf{w}_{ji}-\mathbf{\overline{w}}_i)(\mathbf{w}_{ji}-\mathbf{\overline{w}}_i)'$ which are sufficient statistics for $\mathbf{\Sigma}$, based on the i-th covariate vector. Defining $\mathbf{S}_\mathbf{w}=\sum_{i=1}^n \mathbf{S}_{\mathbf{w}i}$, we have $(\mathbf{\overline{w}}_1,...,\mathbf{\overline{w}}_n,\mathbf{S}_\mathbf{w})$ as the joint sufficient statistics for $(\mathbf{B,\Sigma})$. 
	Conditionally on $\mathbf{\tilde{B}}$ and $\mathbf{\tilde{\Sigma}}$, we have $\mathbf{\overline{w}}_i\sim N_m(\mathbf{\tilde{B}'x}_i,\frac{1}{M}\mathbf{\tilde{\Sigma}})$ and $\mathbf{S}_{\mathbf{w}i}\sim W_m(\mathbf{\tilde{\Sigma}},M-1)$.

	From the $M$ released synthetic data matrices $\mathbf{W}_j\; (j=1,...,M)$, we may define $\mathbf{\overline{W}}_M=\frac{1}{M}\sum_{j=1}^M\mathbf{W}_j$ and define for $\mathbf{B}$ its estimator 
	\begin{equation}\label{eq:par2nd1}
	\mathbf{\overline{B}}^{\bullet}_M=(\mathbf{XX')^{-1}X\overline{W}}_M',
	\end{equation}
	and for $\mathbf{\Sigma}$ its estimator
	\begin{equation}\label{eq:par2nd2}
	\mathbf{S}^{\bullet}_{comb}=\frac{\mathbf{S}_\mathbf{w}+M\times\mathbf{S}^\bullet_{mean}}{Mn-p},
	\end{equation}
	where we define ${\mathbf{S}^\bullet_{mean}}= (\mathbf{\overline{W}}_M-\mathbf{\overline{B}}^{\bullet'}_M\mathbf{X})(\mathbf{\overline{W}}_M-\mathbf{\overline{B}}^{\bullet'}_M\mathbf{X})'$.

		In fact, if the $M$ synthetic datasets are treated as a single synthetic dataset of size $nM$, the estimators obtained for $\mathbf B$ and $\mathbf\Sigma$ 
		will be exactly the same as the ones obtained in (\ref{eq:par2nd1}) and (\ref{eq:par2nd2}). The proof of this fact may be analyzed in Appendix C.

	Analogous to what was done in the previous subsection, one can derive the following inferential results, for $p\geq m$ and $n+\alpha>p+2m+2$.
	
	\begin{enumerate}
		
		\item An UE of $\mathbf{\Sigma}$ will be $\mathbf{\hat{S}}_M=\frac{n+\alpha-p-2m-2}{n-p}\mathbf{S}^{\bullet}_{comb}$ (see Corollary \ref{cor:pdf} Appendix \ref{Aapp:last5}), and for $\alpha=2m+2$, $\mathbf{S}^{\bullet}_{comb}$ will also be an UE for $\mathbf{\Sigma}$.
		
		\item In Corollary \ref{cor:pdf} (see below), we prove that
\vspace{.2cm}
		\begin{equation}\label{eq:T2nd}
		T^\bullet_{comb}=\frac{|(\overline{\mathbf{B}}^{\bullet}_M-\mathbf{B})'(XX')(\mathbf{\overline{B}}^{\bullet}_M-\mathbf{B)}|}{|(Mn-p)\mathbf{S}^{\bullet}_{comb}|}
\vspace{.2cm}
		\end{equation}
		is a pivotal quantity, and that for $\mathbf{A}_1\sim W_m(\mathbf{I}_m,n+\alpha-p-m-1)$, $\mathbf{A}_2\sim W_m(\mathbf{I}_m,n-p)$ and $F_i\sim F_{p-i+1,Mn-p-i+1}\; (i=1,...,m)$, all independent random variables,  
\vspace{.2cm}		$$T^\bullet_{comb}|_{\mathbf{\Omega}}\stackrel{st}{\sim}\left\{\prod_{i=1}^{m}\frac{p-i+1}{Mn-p-i+1}F_{i}\right\}\left|\frac{M+1}{M}\mathbf{I}_{m}+\mathbf{\Omega}\right|,
\vspace{.1cm}
$$
		where $\mathbf{\Omega}$ has the same distribution as $\mathbf{A}_1^{\frac{1}{2}}\mathbf{A}_2^{-1}\mathbf{A}_1^{\frac{1}{2}}$.

		\item If one wants to test a linear combination of the parameters in $\mathbf{B}$, namely, $\mathbf{C = A B }$ where $\mathbf{A}$ is a $k\times p$ matrix with $rank(\mathbf{A})=k\leq p$ and $k\geq m$, one may define 
\vspace{.2cm}
$$T^\bullet_{comb,\mathbf{C}} =\linebreak {\frac{|(\mathbf{A \overline{B}}^{\bullet}_M-\mathbf{C})'(\mathbf{A(XX')^{-1}A'})^{-1}(\mathbf{A\overline{B}}^{\bullet}_M-\mathbf{C})|}{|(Mn-p)\mathbf{\overline{S}}^{\bullet}_{comb}|}},
\vspace{.2cm}
$$
and proceed by noting that
\vspace{.2cm}
		\begin{equation}\label{eq:TC2}
		T^\bullet_{comb,\mathbf{C}}|\mathbf{W}\stackrel{st}{\sim} \left\{\prod_{i=1}^{m}\frac{k-i+1}{Mn-p-i+1}F_{k,i}\right\}\Bigl|\frac{M+1}{M}\mathbf{I}_{m}+\mathbf{\Omega}\Bigr|,
\vspace{.2cm}
		\end{equation}
		with $F_{k,i}\sim F_{k-i+1,Mn-p-i+1}$ being independent random variables and $\mathbf{\Omega}$ defined as in the previous item.
		
		(i)\textit{Test for the significance of $\mathbf{C}$:} in order to test $H_0:\mathbf{C}=\mathbf{C}_0$ versus 
		$H_1:\mathbf{C}\neq\mathbf{C}_0$, we reject $H_0$ whenever $T^\bullet_{comb,\mathbf{C}_0}$ exceeds $\delta_{comb,k,m,p,n;\gamma}$ where 
		$\delta_{comb,k,m,p,n;\gamma}$ satisfies $(1-\gamma)= Pr(T^\bullet_{comb,\mathbf{C}_0}\leq\delta_{comb,k,m,p,n;\gamma})$ when $H_0$ is true. To perform a test for $\mathbf{B=B}_0$ one has to take $\mathbf{A=I}_p$.
		
		(ii)\textit{Confidence set for $\mathbf{C}$:} a $(1-\gamma)$ level confidence set for $\mathbf{C}$ is given by
		\begin{equation}\label{eq:elipsoidAB2}
		\Delta_{comb}(\mathbf{C})=\lbrace\mathbf{C}:T^\bullet_{comb,\mathbf{C}}\leq\delta_{comb,k,m,n,p;\gamma}\rbrace,
		\end{equation}
		where the value of $\delta_{comb,k,m,n,p;\gamma}$ can be obtained by simulating the distribution in (\ref{eq:TC2}).
	\end{enumerate}
	
	Results in 1-3 are derived based on the following Corollaries \ref{cor:pdf} and \ref{cor:dist}, of Theorems \ref{thm:pdf} and \ref{thm:dist}, respectively.
	
	\begin{cor}\label{cor:pdf}
		
		The joint pdf of $\mathbf{\overline{B}}^{\bullet}_M,\mathbf{S}_{comb}^\bullet$ and $\mathbf{\tilde{\Sigma}^{-1}}$, for $\mathbf{\overline{B}}^{\bullet}_M$ and $\mathbf{S}_{comb}^\bullet$ defined in (\ref{eq:par2nd1}) and (\ref{eq:par2nd2}), is proportional to
		\[\begin{array}{l}
		e^{-\frac{1}{2}tr\lbrace\mathbf{(\frac{M+1}{M}\tilde{\Sigma}+\Sigma)}^{-1}(\mathbf{\overline{B}}^{\bullet}_M-\mathbf{B})'\mathbf{XX'}(\mathbf{\overline{B}}^{\bullet}_M-\mathbf{B})+(Mn-p)\mathbf{\tilde{\Sigma}^{-1}}\mathbf{S}^{\bullet}_{comb}\rbrace}\\ 
		~~~~~~~\times \frac{\vert\mathbf{S}^{\bullet}_{comb}\vert^{\frac{Mn-p-m-1}{2}}}{\vert\mathbf{\tilde{\Sigma}}\vert^{\frac{Mn-p+n+\alpha}{2}-m-1}}\vert\mathbf{\Sigma}\vert^{-\frac{n}{2}} |\frac{M}{M+1}\mathbf{\tilde{\Sigma}^{-1}}+\mathbf{\Sigma^{-1}}|^{-p/2} |\mathbf{\tilde{\Sigma}^{-1}+\Sigma^{-1}}|^{-\frac{2n+\alpha-2p-m-1}{2}}\,,
		\end{array}\]
		so that $\overline{\mathbf{B}}^\bullet_M$ and $\mathbf{S}^\bullet_{comb}$, given $\tilde{\mathbf{\Sigma}}$, are independent, with $$\mathbf{\overline{B}}^{\bullet}_M|_{\tilde{\mathbf{\Sigma}}}\sim N_{pm}\left(\mathbf{B},\left(\frac{M+1}{M}\mathbf{\tilde{\Sigma}+\Sigma}\right)\otimes(\mathbf{XX}')^{-1}\right)$$
		and 
		$$\mathbf{S}^{\bullet}_{comb}|_{\tilde{\mathbf{\Sigma}}}\sim W_m\left(\frac{1}{Mn-p}\mathbf{\tilde{\Sigma}},M(n-p)\right).$$
	\end{cor}
	
	\begin{proof}
		See Appendix \ref{App:A}.
	\end{proof}
	
	\begin{cor}\label{cor:dist}
		The distribution of the statistic $T^\bullet_{comb}$ defined in (\ref{eq:T2nd}) can be obtained from the decomposition
		$$T^\bullet_{comb}|_\mathbf{\Omega}\stackrel{st}{\sim}\left\{\prod_{i=1}^{m}\frac{p-i+1}{Mn-p-i+1}F_{i}\right\}\left|\frac{M+1}{M}\mathbf{I}_{m}+\mathbf{\Omega}\right|$$
		where $F_i\sim F_{p-i+1,Mn-p-i+1}$ are independent random variables, themselves independent of $\mathbf{\Omega}$, which has the same distribution as $\mathbf{A}_1^{\frac{1}{2}}\mathbf{A}_2^{-1}\mathbf{A}_1^{\frac{1}{2}}$ with $\mathbf{A}_1\sim W_m(\mathbf{I}_m, n+\alpha-p-m-1)$ and $\mathbf{A}_2\sim W_m(\mathbf{I}_m, n-p)$, two independent random variables.
	\end{cor}
	
	\begin{proof}
		See Appendix \ref{App:A}.
	\end{proof}

		\begin{rem}
			Similar to what happens with the statistic $T^\bullet_M$ in (\ref{eq:T1st}), the statistic $T^\bullet_{comb}$ in (\ref{eq:T2nd}) also degenerates towards zero when $n\rightarrow\infty$ or $M\rightarrow\infty$, and similarly to what happens with $T^\bullet_M$,
			$$
			(Mn-p)^m\,T^\bullet_{comb}|_{\mathbf{\Omega}}\xrightarrow[n\rightarrow\infty]{d}\left\{\prod_{i=1}^{m}\chi^2_{p-i+1}\right\}\left|\frac{M+1}{M}\mathbf{I}_{m}+\mathbf{\Omega}\right|
			$$
			and
			$$
			(Mn-p)^m\,T^\bullet_{comb}|_{\mathbf{\Omega}}\xrightarrow[M\rightarrow\infty]{d}\left\{\prod_{i=1}^{m}\chi^2_{p-i+1}\right\}\left|\mathbf{I}_{m}+\mathbf{\Omega}\right|.
			$$
			Using the simple scale change $T^\bullet_{comb2}=(Mn-p)^m\,T^\bullet_{comb}=\frac{|(\overline{\mathbf{B}}^{\bullet}_M-\mathbf{B})'(XX')(\mathbf{\overline{B}}^{\bullet}_M-\mathbf{B)}|}{|\mathbf{\overline{S}}^{\bullet}_{comb}|}$
			one would have
			$$
			T^\bullet_{comb2}|_{\mathbf{\Omega}}\xrightarrow[n\rightarrow\infty]{d}\left\{\prod_{i=1}^{m}\chi^2_{p-i+1}\right\}\left|\frac{M+1}{M}\mathbf{I}_{m}+\mathbf{\Omega}\right|
			$$
			and
			$$
			T^\bullet_{comb2}|_{\mathbf{\Omega}}\xrightarrow[M\rightarrow\infty]{d}\left\{\prod_{i=1}^{m}\chi^2_{p-i+1}\right\}\left|\mathbf{I}_{m}+\mathbf{\Omega}\right|,
			$$
			similar to what happens with $T^\bullet_M$.
		\end{rem}

\section{SIMULATION STUDIES}
\label{sec:sim}
	
	In order to compare the PPS and the FPPS methods with the Plug-in Sampling method we present the results of some simulations analogous to the ones presented in \cite{moura16}. The objectives of these simulations are: (i) to show that the inference methods developed in Section \ref{sec:post} perform as predicted, and (ii) to compare the measures (\textit{radius}) obtained from our methods with the ones from the Plug-in method. All simulations were carried out using the software Mathematica$^\circledR$. To conduct the simulation, we take the population distribution as a multivariate normal distribution with expected value given by the right hand side of (\ref{eq:model}), for $m=2$ and $p=3$, with matrix of regressor coefficients
	$$\mathbf{B}=\left(
	\begin{matrix}
	1 & 2 \\
	3 & 2 \\
	1 & 1
	\end{matrix}
	\right)$$
	and covariance matrix
	$$\mathbf{\Sigma}=\left(
	\begin{matrix}
	1 & 0.5 \\
	0.5 & 1 \\ \end{matrix} \right).$$ 
	We set $\alpha=6$ in order to have both $\mathbf{\bar{S}}^\bullet_M$ and $\mathbf{S}^\bullet_{comb}$ as the unbiased estimators of $\mathbf{\Sigma}$. The regressor variables $x_{1i},x_{2i}, x_{3i}, i=1,...,n$ are generated as i.i.d. $N(1,1)$ 
	and held fixed for the entire simulation. Based on Monte Carlo simulation with $10^5$ iterations, 
	we compute an estimate of the coverage probability of the confidence regions for $\mathbf{B}$ and $\mathbf{C=AB}$ given by (\ref{eq:elipsoidAB}) and (\ref{eq:elipsoidAB2}), defined as percentage of observed values of the statistics smaller than the respective theoretical cut-off points, with $\mathbf{A}=\left(\begin{smallmatrix} 0 & 1 & 0\\ 0 & 0 & 1 \end{smallmatrix}\right)$, using the methodologies described in Section \ref{sec:post}. For $M=1$, $M=2$ and $M=5$, the estimated coverage probabilities of the confidence sets are shown in Table 2 under the columns $\mathbf{B}(1)$ and $\mathbf{AB}(1)$ for the first new procedure in Subsection \ref{ssec:mul1}, and under the columns $\mathbf{B}(2)$ and $\mathbf{AB}(2)$ for the second new procedure in Subsection \ref{ssec:mul2}. For $M=1$, a single column is shown for each confidence region since the two new procedures are the same.
	
	\begin{table}[h!]
		\caption{\footnotesize{Average coverage for $\mathbf{B}$ and $\mathbf{AB}$}}
		\footnotesize
		\vspace{-10pt}
		\begin{center}
			\footnotesize{
				\begin{tabular}{| r | c | c | c | c | c | c | c | c | c | c | }\hline
					\multirow{3}{*}{$n$} & \multicolumn{2}{c|}{$M=1$} & \multicolumn{4}{c|}{$M=2$} & \multicolumn{4}{c|}{$M=5$} \\ \cline{2-11}
					&  \multirow{2}{*}{$\mathbf{B}$} & \multirow{2}{*}{$\mathbf{AB}$}  & \multicolumn{2}{c|}{1st Approach} & \multicolumn{2}{c|}{2nd Approach} & \multicolumn{2}{c|}{1st Approach} & \multicolumn{2}{c|}{2nd Approach} \\
					&  &  & $\mathbf{B} (1)$ & $\mathbf{AB} (1)$ & $\mathbf{B} (2)$ & $\mathbf{AB} (2)$ & $\mathbf{B} (1)$ & $\mathbf{AB} (1)$ & $\mathbf{B} (2)$ & $\mathbf{AB} (2)$\\ \hline
					10  & 0.949  & 0.951    & 0.949  & 0.949 & 0.951 & 0.949 & 0.951 & 0.950 & 0.949 & 0.951 \\[-1pt]
					50  & 0.949  & 0.950    & 0.951  & 0.951 & 0.950 & 0.951 & 0.951 & 0.950 & 0.949 & 0.948  \\[-1pt]
					100 & 0.949  & 0.949    & 0.951  & 0.950 & 0.949 & 0.951 & 0.949 & 0.951 & 0.951 & 0.950  \\[-1pt]
					200 & 0.951  & 0.951    & 0.949  & 0.951 & 0.951 & 0.949 & 0.950 & 0.951 & 0.950 & 0.951  \\ \hline
				\end{tabular}}
			\end{center}
			\vspace{-20pt}
		\end{table}
		
		The results reported in Table 2 for samples of size $n=10, 50, 100, 200$, show that, based on singly and multiply imputed synthetic data, the 0.95 confidence sets for 
		$\mathbf{B}$ and $\mathbf{AB}$ have an estimated coverage probability approximately equal to 0.95, confirming that the confidence sets perform as predicted.
		
		In order to measure the {\it radius} (distance between the center and the edge) of the confidence sets, we use the same measure proposed in \cite{moura16}, which is
		$$\Upsilon_M=d^*_{M,m,n,p,\gamma}\times |\tilde{\mathbf{S}}^\bullet_M|,$$ 
		where $d^*_{M,m,n,p,\gamma}$ is the cut-off point in (\ref{eq:elipsoidAB}) or (\ref{eq:elipsoidAB2}). Here we take $M=0$ for the original data, with
		$\tilde{\mathbf{S}}^\bullet_0=(n-p)\mathbf{S}$, $M=1$ for the singly imputed synthetic data and $M=2,5$ for the multiply imputed synthetic data, with  
		$\tilde{\mathbf{S}}^\bullet_M=M(n-p)\mathbf{\overline{S}}^{\bullet}_M$ for the first new procedure, 
		and $\tilde{\mathbf{S}}^\bullet_M=(Mn-p)\mathbf{S}^\bullet_{comb}$ for the second new procedure. 
		The expected value of this measure will be
		$$E(\Upsilon_M)=d^*_{M,m,n,p,\gamma}\times\frac{(n-p)!}{(n-p-m)!}\times K_{M,n,p,m} |\mathbf{\Sigma}|$$
		where $K_{0,n,p,m}=1$ for the original data, $$K_{M,n,p,m}=\frac{(-2+\kappa_{n,p,\alpha,m}-m)!}{(-2+\kappa_{n,p,\alpha,m})!}\frac{(Mn-Mp)!}{(Mn-Mp-m)!}$$ 
		for the procedure in Subsection \ref{ssec:mul1} and $$K_{M,n,p,m}=\frac{(-2+\kappa_{n,p,\alpha,m}-m)!}{(-2+\kappa_{n,p,\alpha,m})!}\frac{(Mn-p)!}{(Mn-p-m)!}$$ 
		for the procedure in Subsection \ref{ssec:mul2}, where $\kappa_{n,\alpha,p,m}=n+\alpha-p-m-1$, assuming $n+\alpha>p+2m+2$.
		For more details about these expected values we refer to Appendix \ref{Aapp:last6}.

		We present in Table 3 the average of the simulated values of the \textit{radius} $\Upsilon_M$ and its expected value $E(\Upsilon_M)$ for the confidence sets $\Delta_M(\mathbf{B})$ (first procedure) and $\Delta_{comb}(\mathbf{B})$ (second procedure), and in Table 4 the same values for the confidence sets $\Delta_M(\mathbf{C})$ (first procedure) and $\Delta_{comb}(\mathbf{C})$ (second procedure), for 
		$M=0,1,2,5$ and $n=10,50,200$. These values may be compared with the values obtained in \cite{moura16} for the Plug-in Sampling.

			\begin{table}[h!]
				\caption{Average values of $\Upsilon_M$ and the values of $E(\Upsilon_M)$ for the confidence set for $\mathbf{B}$.}
				\footnotesize
				\vspace{-10pt}
				\begin{center}
				
				\begin{tabular}{c}
					\begin{tabular}{| r | c || c | c || c | c | c | c | }\hline
						\multirow{3}{*}{$n$} &	\multirow{3}{*}{Orig} & \multicolumn{2}{c||}{$M=1$} & \multicolumn{4}{c|}{$M=2$}   \\\cline{3-8}
						& & \multirow{2}{*}{avg} & \multirow{2}{*}{exp} & \multicolumn{2}{c|}{1st Procedure} & \multicolumn{2}{c|}{2nd Procedure} \\
						&  & &     & avg     & exp     & avg     & exp    \\ \hline
						10  & 36.97  & 507.25 & 512.19 & 251.55 & 252.55 & 237.64 & 238.68   \\
						50  & 19.11  & 176.36 & 176.53 & 121.23 & 121.52 & 121.23 & 121.48  \\
						200 & 17.52  & 154.93 & 156.06 & 105.81 & 106.61 & 105.90 & 106.72  \\ \hline
					\end{tabular}
				\\
				\vspace{-5pt}
				\\
					\begin{tabular}{| r | r | r | r | r | }\hline
						\multirow{3}{*}{$n$}  & \multicolumn{4}{c|}{$M=5$}  \\\cline{2-5}
						& \multicolumn{2}{c|}{1st Procedure} & \multicolumn{2}{c|}{2nd Procedure} \\
						&   avg     & exp     & avg     & exp     \\ \hline
						10  &175.34 & 176.18 & 163.82 & 168.92  \\
						50  &92.25 & 92.80 & 92.28 & 92.84 \\
						200 & 81.89 & 82.39 & 81.91 & 82.40 \\ \hline
					\end{tabular}
				
			\end{tabular}
			\end{center}
				\vspace{-10pt}
			\end{table}

			\begin{table}[h!]
			\caption{Average values of $\Upsilon_M$ and the values of $E(\Upsilon_M)$ for the confidence set for $\mathbf{C=AB}$.}
			\footnotesize
			\vspace{-10pt}
				\begin{center}
					
					\begin{tabular}{c}
						\begin{tabular}{| r | r || r | r || r | r | r | r | }\hline
							\multirow{3}{*}{$n$} &	\multirow{3}{*}{Orig} & \multicolumn{2}{c||}{$M=1$} & \multicolumn{4}{c|}{$M=2$}   \\\cline{3-8}
							& & \multirow{2}{*}{avg} & \multirow{2}{*}{exp} & \multicolumn{2}{c|}{1st Procedure} & \multicolumn{2}{c|}{2nd Procedure} \\
							&  & &     & avg     & exp     & avg     & exp    \\ \hline
							10  & 13.43 & 172.64 & 172.32 & 92.23 & 92.44 & 86.24 & 86.61   \\
							50  & 7.33 & 68.93 & 68.99 & 47.75 & 47.86 & 47.45 & 47.55  \\
							200 & 7.10 & 60.65 & 61.09 & 41.74 & 42.05 & 41.74 & 42.05  \\ \hline
						\end{tabular}
						\\
						\vspace{-5pt}
						\\
						\begin{tabular}{| r | c | c | c | c | }\hline
							\multirow{3}{*}{$n$}  & \multicolumn{4}{c|}{$M=5$}  \\\cline{2-5}
							& \multicolumn{2}{c|}{1st Procedure} & \multicolumn{2}{c|}{2nd Procedure} \\
							&   avg     & exp     & avg     & exp     \\ \hline
							10  & 63.07 & 63.38 & 61.34 & 61.74  \\
							50  & 35.32 & 35.52 & 35.08 & 35.27 \\
							200 & 32.47 & 32.51 & 32.54 & 32.53 \\ \hline
						\end{tabular}
						
					\end{tabular}
				\end{center}
				\vspace{-10pt}
			\end{table}

		Observing Tables 3 and 4 and comparing the entries in these tables with the results in \cite{moura16} for Plug-in Sampling, we may see that when synthetic data are generated under FPPS, larger \textit{radius} are obtained. In the singly imputed case, one can observe that the PPS synthetic datasets will lead to a \textit{radius} that is approximately two and half times that of the \textit{radius} under Plug-in Sampling.  As the number $M$ of released synthetic datasets increases, $\Upsilon_M$ slowly decreases, increasing however the difference of the \textit{radius} between the FPPS and the Plug-in methods.
		Eventually, one may need very large values of $M$, in order to have values of $\Upsilon_M$ close to the value of $\Upsilon_0$. As in \cite{moura16} we also observe that the values of $\Upsilon_M\; (M>1)$, for both new FPPS procedures become identical for larger sample sizes.

\section{AN APPLICATION USING CURRENT POPULATION\\ SURVEY DATA}
\label{sec:app}

		In this section, we provide an application based on the same real data used in \cite{moura16} to compare the original data inference with the one obtained via PPS, for the single imputation case, and via FPPS, for the multiple imputation case. The data are from the U.S. 2000 Current Population Survey (CPS) March supplement, available online at http://www.census.gov.cps/. Further details on the data may be found in \cite{moura16}.
		
		In this application, $\mathbf{x}$, the vector of regressor variables, is defined as
		\[
\begin{array}{l}
		\mathbf{x}= \Big(1,N,L,A,I(E=34),...,I(E=37),I(E=39),...,I(E=46), \medskip \\ 
		\hskip1.5cm I(M=3),...,I(M=7),I(R=2),I(R=4),I(S=2)\Big)', 
		\end{array}
		\]
		where N, L, A,  are respectively, the number of people in household, the number of people in the household who are less than 18 years old and the age for the head of household, E, M, R and S, are respectively, the education level for the head of the household (coded to take values 31, 34-37, 39-46), the marital status for the head of the household (coded to take values 1,3-7), the race of the head of the household (coded to take values 1,2,4) and the sex of the head of the household (coded to take values 1,2). $\rm I(E=34)$ is the indicator variable for $E=34$, $\rm I(E=35)$ is the indicator variable for $E=35$, and so on, and where the indicator variable for the first code present in the sample for each variable is taken out in order to make the model matrix full rank.
		The vector $\mathbf{y}$ of response variables will be formed by the same three numerical variables used in \cite{moura16}, namely, \textit{total household income}, \textit{household alimony payment} and \textit{household property tax}. After deleting all entries where at least one of these variables are reported as 0, we were left with a sample size of 141, and as such the model matrix $\mathbf{X}=[\mathbf{x}_1\cdots \mathbf{x}_n]$ has thus $p=24$ rows, $n=141$ columns, with rank equal to 24. Throughout this section we will assume $\alpha=8$ in order to have $\mathbf{S}^\bullet_M$ and $\mathbf{S}^\bullet_{comb}$ as unbiased estimators of $\mathbf{\Sigma}$. Via PPS method we generate a single synthetic dataset and show in expression (\ref{exp:S}) the realizations of the unbiased estimator $\mathbf{S}^\bullet$ for $\mathbf{\Sigma}$ and of the estimator $\mathbf{S}$ for the original data, respectively denoted by $\widetilde{\mathbf{S}}^\bullet_1$ and $\mathbf{\widetilde{S}}$
		\begin{equation}\footnotesize\label{exp:S}
		\mathbf{\widetilde{S}^{\bullet}}_1=\left(\!
		\begin{array}{rrr}
		1.58572 & -0.20443 & 0.27981 \\
		-0.20443 & 1.61395 & 0.16089 \\
		0.27981 & 0.16089 & 0.34648 \\
		\end{array}
		\right),
		\;\;
		\mathbf{\widetilde{S}}=\left(\!
		\begin{array}{rrr}
		1.1980 & -0.0375 & 0.2970 \\[-2pt]
		-0.0375 & 1.0699 & 0.1175 \\[-2pt]
		0.2970 & 0.1175 & 0.4045 \\
		\end{array}
		\right).
		\end{equation}
		In Table 5 we show the realizations of the unbiased estimator $\mathbf{B}^\bullet_1$ of $\mathbf{B}$ and of the estimator $\hat{\mathbf{B}}$ of the original data, respectively denoted by $\widetilde{\mathbf{B}}^\bullet_1$ and $\mathbf{\widetilde{\hat{B}}}$.
		
		\begin{table}[h!]
			
			\caption{Estimates of the regressor coefficients from the FPPS synthetic data ($\mathbf{\widetilde{B}}^{\bullet}$), Plug-in synthetic data ($\mathbf{\widetilde{B}}^*$)  and from the original data.}
			\scriptsize
			\vspace{-0.4cm}
			\begin{center}$\begin{array}{c||rrr||rrr|||rrr}
				& \multicolumn{3}{c||}{FPPS} &\multicolumn{3}{c|||}{Plug-in} &  \\ [-3pt]
				\multirow{2}{*}{regressor} & \multicolumn{3}{c||}{Synthetic Data\; (\mathbf{\widetilde{B}}^{\bullet})} &\multicolumn{3}{c|||}{Synthetic Data\; (\mathbf{\widetilde{B}}^{*})} & \multicolumn{3}{c}{Original Data \;(\mathbf{\widetilde{\hat{B}}})} \\
				& \multicolumn{1}{c}{\rm I} & \multicolumn{1}{c}{\rm AP} & \multicolumn{1}{c||}{\rm PT} & \multicolumn{1}{c}{\rm I} & \multicolumn{1}{c}{\rm AP} & \multicolumn{1}{c|||}{\rm PT} & \multicolumn{1}{c}{\rm I} & \multicolumn{1}{c}{\rm AP} & \multicolumn{1}{c}{\rm PT}\\ \hline
				\text{Intercept} & 11.4996 & 3.3381 & 8.1713 & 10.1829 & 3.7094 & 10.9787 & 9.8339 & 4.6663 & 10.1095 \\[0pt]
				\text{N} & 0.2801 & -0.2562 & 0.6317 &-0.0938 & 0.1435 & 0.6189 & 0.0457 & 0.0375 & 0.4585 \\[0pt]
				\text{L} &-0.3996 & 0.4960 & -0.6017 & 0.0812 & 0.0163 & -0.5932 & 0.0186 & 0.1310 & -0.3851 \\[0pt]
				\text{A} &  -0.0061 & 0.0223 & 0.0018  & 0.0075 & 0.0285 & -0.0097 & 0.0118 & 0.0181 & -0.0020 \\[0pt]
				\text{I(E=34)} & -4.7732 & 0.3476 & -0.4662 & -6.6680 & 1.2055 & -2.0664 & -4.4348 & 0.5944 & -1.2291 \\[0pt]
				\text{I(E=35)} & -5.5990 & 2.8081 & 1.9914 & -1.2231 & -0.0154 & -0.7091 & -1.4060 & 0.9188 & -0.1468 \\ [0pt]
				\text{I(E=36)} & -4.2467 & 2.2712 & 0.6907 & -0.4478 & 2.1718 & -0.9172 & -2.3100 & 1.0416 & -0.5002 \\[0pt]
				\text{I(E=37)} &  -3.5281 & 0.7339 & 1.4653 & -1.1547 & 1.3009 & -1.0659 & -2.0490 & 0.7410 & 0.2335 \\[0pt]
				\text{I(E=39)} & -3.3369 & 1.5590 & 1.0109  & -2.5737 & 0.7234 & -1.1346 & -2.2208 & 0.4054 & -0.4136 \\[0pt]
				\text{I(E=40)} & -2.8766 & 1.7608 & 1.2350 & -1.8032 & 1.0617 & -0.6940 & -1.8834 & 0.8519 & 0.0852 \\[0pt]
				\text{I(E=41)} & -2.8266 & 2.7954 & 2.3165 & -1.5615 & 1.6881 & -0.0291 & -1.9468 & 1.4222 & 0.1094 \\[0pt]
				\text{I(E=42)} & -3.5901 & 2.3990 & 0.7908 & -2.4543 & 2.0378 & -1.1494 & -2.3381 & 1.3840 & -0.0808 \\[0pt]
				\text{I(E=43)} &-1.9852 & 2.1149 & 1.9765& -1.7090 & 1.1722 & -0.4341 & -1.5057 & 1.0766 & 0.5309 \\[0pt]
				\text{I(E=44)} & -3.2012 & 2.0495 & 1.7665 & -2.2668 & 1.5629 & -0.2140 & -1.8082 & 1.1301 & 0.4936 \\[0pt]
				\text{I(E=45)} &   0.1813 & 1.1103 & 1.7535 & -1.8984 & 2.1024 & -0.4636 & -0.9893 & 0.7958 & 0.3057 \\[0pt]
				\text{I(E=46)} & 0.5791 & 2.3091 & 3.5534 &  0.4558 & 1.4836 & 1.1497 & -0.6198 & 1.0766 & 1.0624 \\[0pt]
				\text{I(M=3)} &-2.3691 & 0.8545 & -0.3594 & -1.9077 & -0.4988 & -0.4836 & -2.7258 & 0.0964 & -0.2156 \\[0pt]
				\text{I(M=4)} & -4.4234 & 2.2640 & -1.2282 & -0.0088 & 0.5609 & -0.2349 & -0.0134 & 0.5887 & 0.3864 \\[0pt]
				\text{I(M=5)} & -1.0787 & 1.5611 & 0.1170 & 0.3767 & 0.6729 & 0.1184 & 0.1455 & 0.4770 & 0.1558 \\[0pt]
				\text{I(M=6)} & -0.8300 & -0.2358 & -0.2713 & 0.3948 & -0.3092 & -0.1046 & -0.7122 & -0.4448 & -0.4025 \\[0pt]
				\text{I(M=7)} & -2.8242 & 2.9533 & 0.5456  & 1.0576 & 0.5476 & 0.5187 & -0.1990 & 1.1750 & 0.6685 \\[0pt]
				\text{I(R=2)} & 0.3378 & 3.8443 & 1.4196 & -1.0805 & 3.0078 & -0.1619 & -0.9205 & 1.3432 & 0.4696 \\[0pt]
				\text{I(R=4)} &  0.0340 & 1.9168 & -0.4519 & 0.6883 & -0.3211 & 0.3639 & -0.7040 & 0.0975 & -0.1618 \\[0pt]
				\text{I(S=2)} &1.3582 & -0.4793 & -0.1588  & 0.0564 & -0.2309 & -0.2849 & 0.1236 & -0.1355 & -0.4025 \\[0pt]
				\end{array}$
			\end{center}
			\vspace{-10pt}
		\end{table}
		
		At a first glance the estimates originated via Plug-in Sampling (see \cite{moura16}) seem to be more in agreement with the original data estimates than the ones drawn from PPS. Nevertheless, this is only one draw and it could be a question of chance to originate `better' or `worse' data. Therefore, one must conduct inferences on the regression coefficients based on multiple draws.
		
		Inferences on regression coefficients are obtained by applying the methodologies in Subsections \ref{ssec:mul1} and \ref{ssec:mul2}, to analyze the
		singly imputed synthetic dataset and multiply imputed synthetic datasets, considering $M=1$, $M=2$ and $M=5$, using the statistics $T^\bullet_M$ and $T^\bullet_{comb}$ and their empirical distributions based on simulations with $10^4$ iterations, to test the fit of the model and the significance of some regressors for $\gamma=0.05$. Regarding the test of fit of the model one will find, for all values of $M$, results equivalent to the ones obtained for the case when synthetic data are generated via Plug-in Sampling, i.e., concluding that the explanatory variables in $\mathbf{x}$ have a significant role in determining the values of the response variables in $ \mathbf{y} $ since the obtained p-values, computed as the fraction of values of the empirical distribution of the corresponding statistic that are larger than the computed value of the statistic, were all approximately zero. 	
		The cut-off points obtained from the empirical distributions of $T^\bullet_M$ and $T^\bullet_{comb}$ (respectively associated with the first and second procedures in Subsections \ref{ssec:mul1} and \ref{ssec:mul2}) are approximately equal to $0.50357$, for $M=1$ (where first and second procedures coincide), to $0.03460$ and $0.02569$, for $M=2$, and to $0.00149$ and $0.00094$, for $M=5$. In Figure 1, one can see a histogram associated with the empirical distributions of both $T^\bullet_M$ and $T^\bullet_{comb}$ for $M=1,2$ and $5$ (for $m=3$, $p=24$, $n=141$, $\alpha=8$ and $10^4$ simulation sizes), recalling that for $M=1$ these two statistics are the same.

		\begin{figure}[h!]
			\vspace{0.25cm}
			\centering
			\begin{subfigure}{.33\textwidth}
				\centering
				\includegraphics[width=1\linewidth]{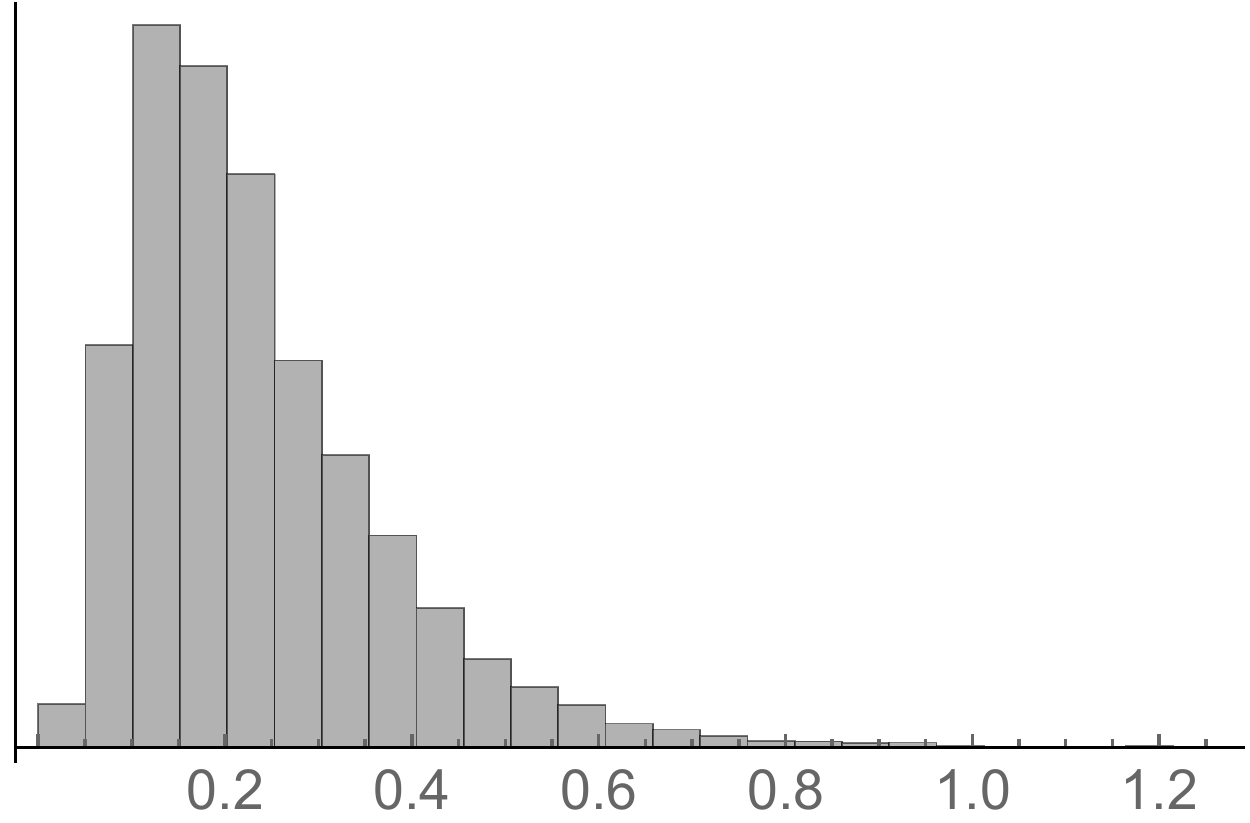}
				\vspace{-15pt}
				\caption{$M=1$}
			\end{subfigure}%
			\begin{subfigure}{.33\textwidth}
				\centering
				\includegraphics[width=1\linewidth]{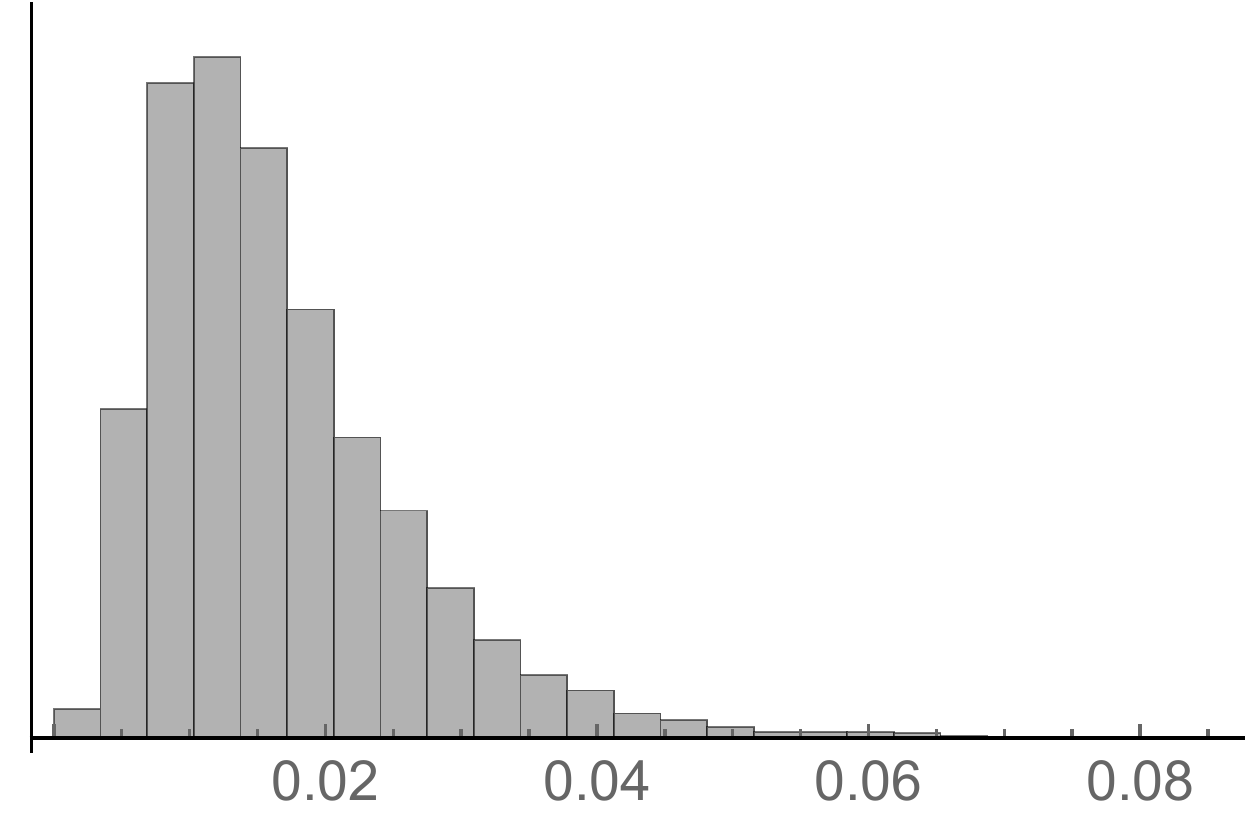}
				\vspace{-15pt}
				\caption{$M=2$ (first procedure)}
			\end{subfigure}
			\begin{subfigure}{.33\textwidth}
				\centering
				\includegraphics[width=1\linewidth]{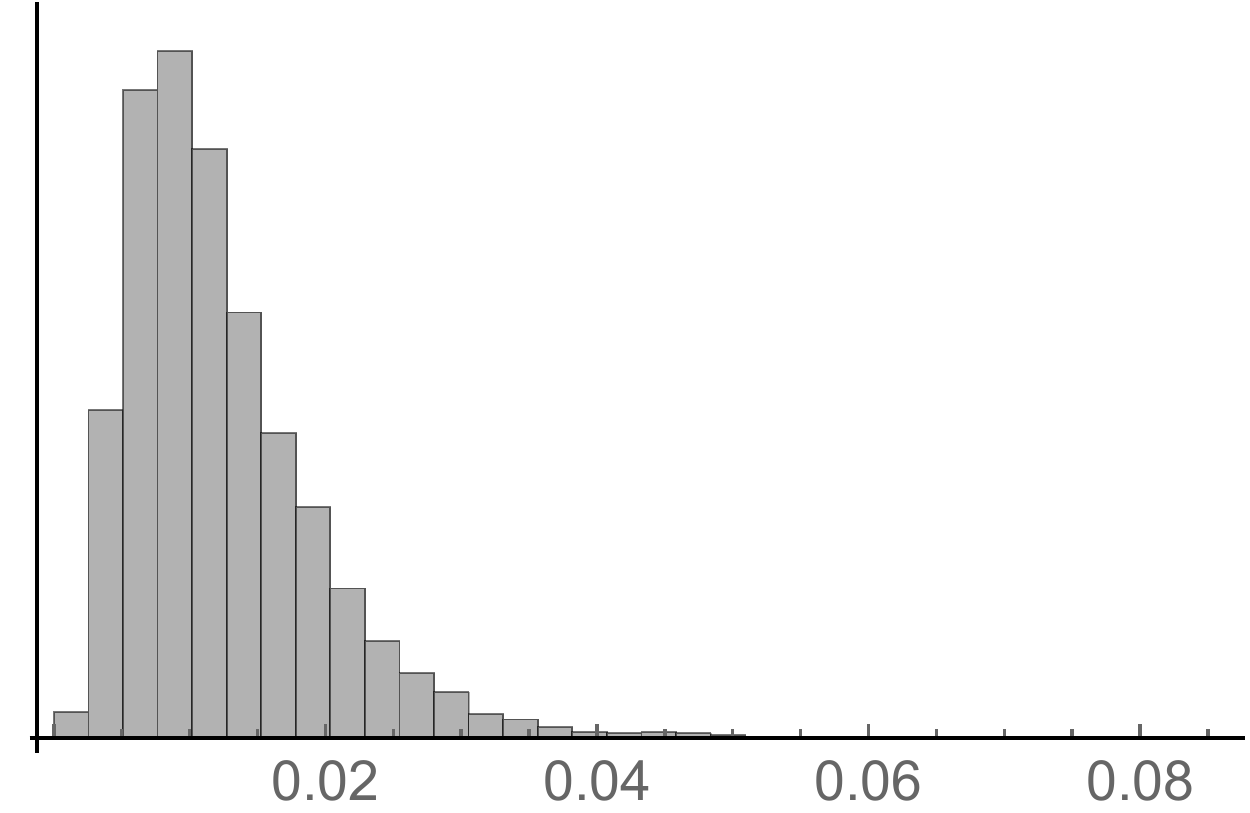}
				\vspace{-15pt}
				\caption{$M=2$ (second procedure)}
			\end{subfigure}
			\begin{subfigure}{.33\textwidth}
				\vspace{10pt}
				\centering
				\includegraphics[width=1\linewidth]{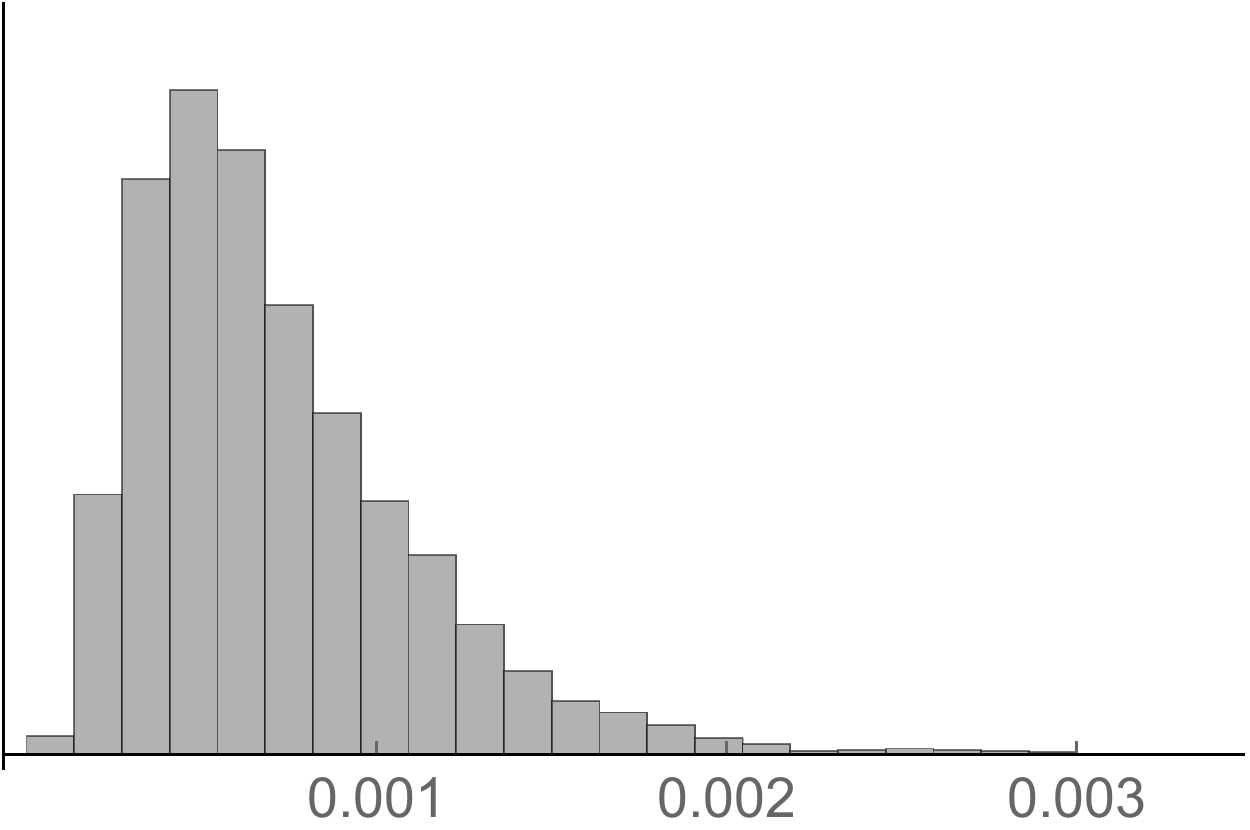}
				\vspace{-15pt}
				\caption{$M=5$ (first procedure)}
			\end{subfigure}
			\begin{subfigure}{.33\textwidth}
				\vspace{10pt}
				\centering
				\includegraphics[width=1\linewidth]{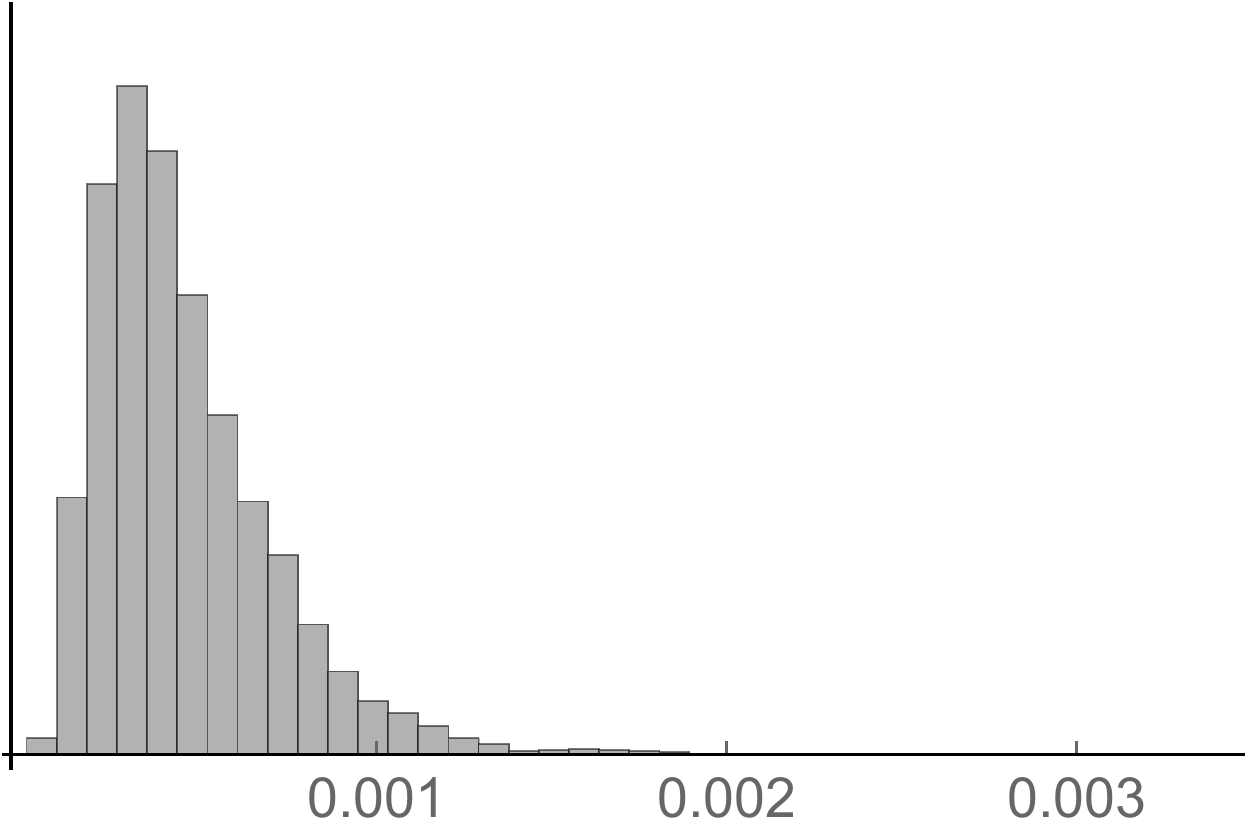}
				\vspace{-15pt}
				\caption{$M=5$ (second procedure)}
			\end{subfigure}
			\vspace{-5pt}
			\caption{
			Histograms (all with same vertical scale) of the empirical distributions of both $T^\bullet_M$ and $T^\bullet_{comb}$ for $M=1,2$ and $5$ (for $m=3$, $p=24$, $n=141$, $\alpha=8$ and $10^4$ simulation sizes)}\label{fig:dist}
		\end{figure}
		
\begin{figure}[h!]
			
			\centering
			\begin{subfigure}{.45\textwidth}
				\centering
				\includegraphics[width=1\linewidth]{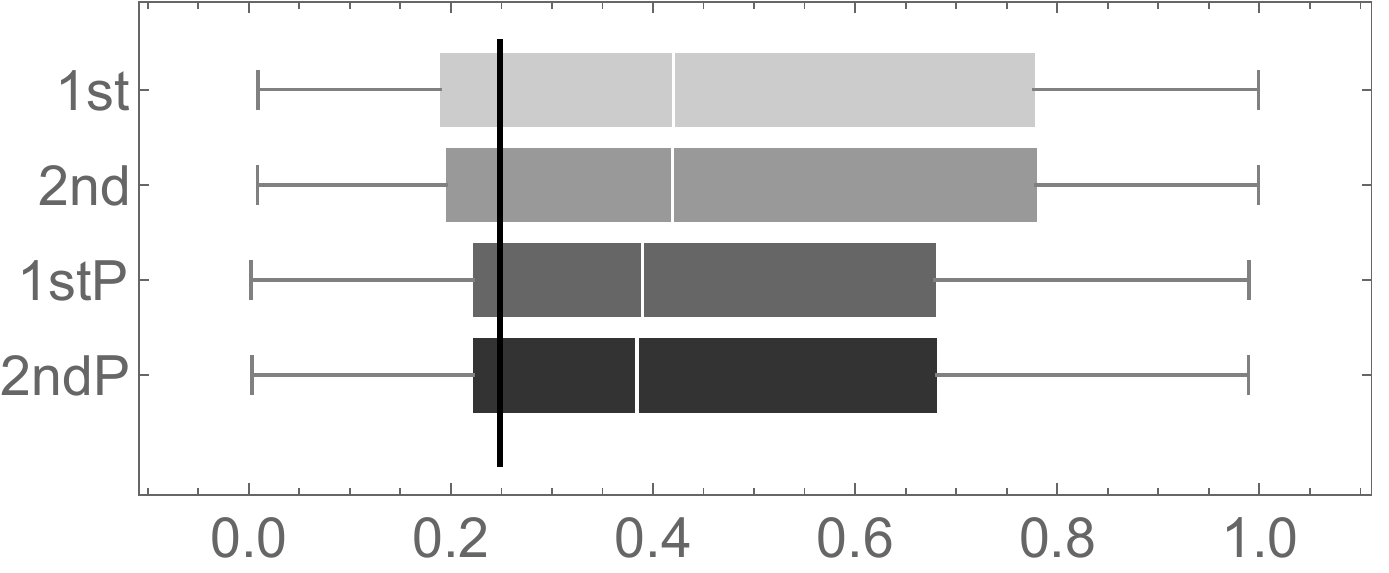}
				\label{fig:3boxM21}
				\vspace{-15pt}
				\caption{$M=2$}
			\end{subfigure}%
			\begin{subfigure}{.45\textwidth}
				\centering
				\includegraphics[width=1\linewidth]{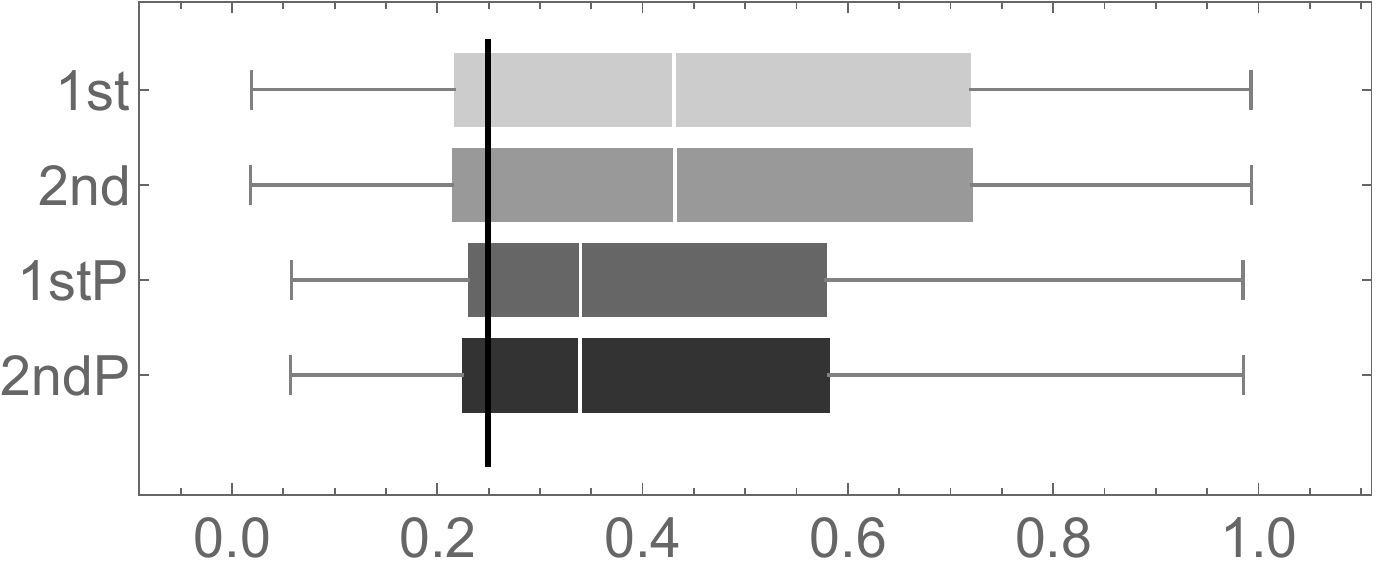}
				\label{fig:3boxM5}
				\vspace{-15pt}
				\caption{$M=5$}
			\end{subfigure}
			
			\begin{subfigure}{.45\textwidth}
				\centering
				\includegraphics[width=1\linewidth]{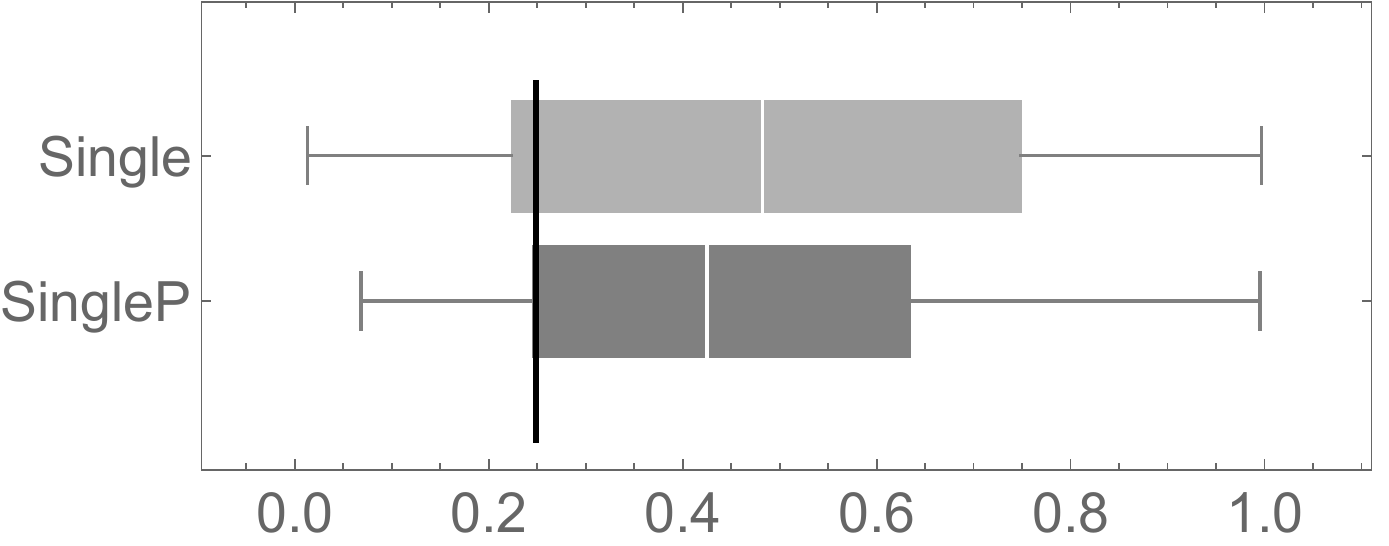}
				\label{fig:3boxSingle}
				\vspace{-15pt}
				\caption{$M=1$}
			\end{subfigure}
			\vspace{-5pt}
			\caption{
			Box-plots of p-values obtained, when testing the joint significance of I(R=2), I(R=4) and I(S=2), from 100 draws of synthetic datasets using procedures in Section \ref{sec:post} and using Plug-in Sampling method from \cite{moura16}, for $M=1$, $M=2$ and $M=5$ .}\label{fig:box1}
		\end{figure}

		In order to test the significance of some regressors, we propose to study two different cases, using in each case the same sets of regressors as in \cite{moura16}. Therefore, we will test the significance of regressor variables R and S, for the first case, and regressor variables A and E, for the second case. As such, in the first case, we will consider a $3\times 24$ matrix 
		$$\mathbf{A}=\left( \begin{array}{c|c}
		\textbf{\Large 0}_{3\times21} & \textbf{\Large I}_3 \end{array} \right)$$
		and we will be interested in testing the hypothesis $H_0: \mathbf{AB=C_0}$, where $\mathbf{C}_0$ is a $3\times 3$ matrix consisting of only zeros. We now generate 100 draws of $M=1$, $M=2$ and $M=5$ synthetic datasets and gather the different p-values obtained when using the statistics in (\ref{eq:T1st}) and (\ref{eq:T2nd}). In Figure \ref{fig:box1}, one may analyze the box-plots of the p-values obtained for each procedure together with the ones obtained in \cite{moura16} for the same sets of variables, where under Single, 1st and 2nd, one has the box-plots associated with the new procedures developed in this paper and under SingleP, 1stP and 2ndP, the box-plots associated to the Plug-in Sampling method. 
		The existing line in the box-plots marks the original data p-value 0.249, obtained using the $T_{O,\mathbf{C}}$ statistic in (3) of \cite{moura16}. 
		It is important to note that in the case of single imputation ($M=1$) the FPPS method reduces to the usual PPS method.

		In general, from Figure \ref{fig:box1}, we may note in both new procedures a larger spread of the p-values when compared with the p-values gathered from Plug-in Sampling, presenting a distribution of p-values with larger values than the original, nonetheless with the majority of these p-values leading to similar conclusions as those obtained from the original data for $\gamma=0.05$, that is, to not reject the null hypothesis that variables R and S do not have significant influence on the response variables.
		
		We may note that in general, in cases where the p-value obtained from the original data is rather low, we expect to obtain larger p-values for the synthetic data, given the inherent variability of these synthetic data and the ``need'' of the inferential exact methods to preserve the $1-\gamma$ coverage level, and impossibility of compressing the synthetic data p-values towards zero.
		
\begin{figure}[h!]
			
			\centering
			\begin{subfigure}{.45\textwidth}
				\centering
				\includegraphics[width=1\linewidth]{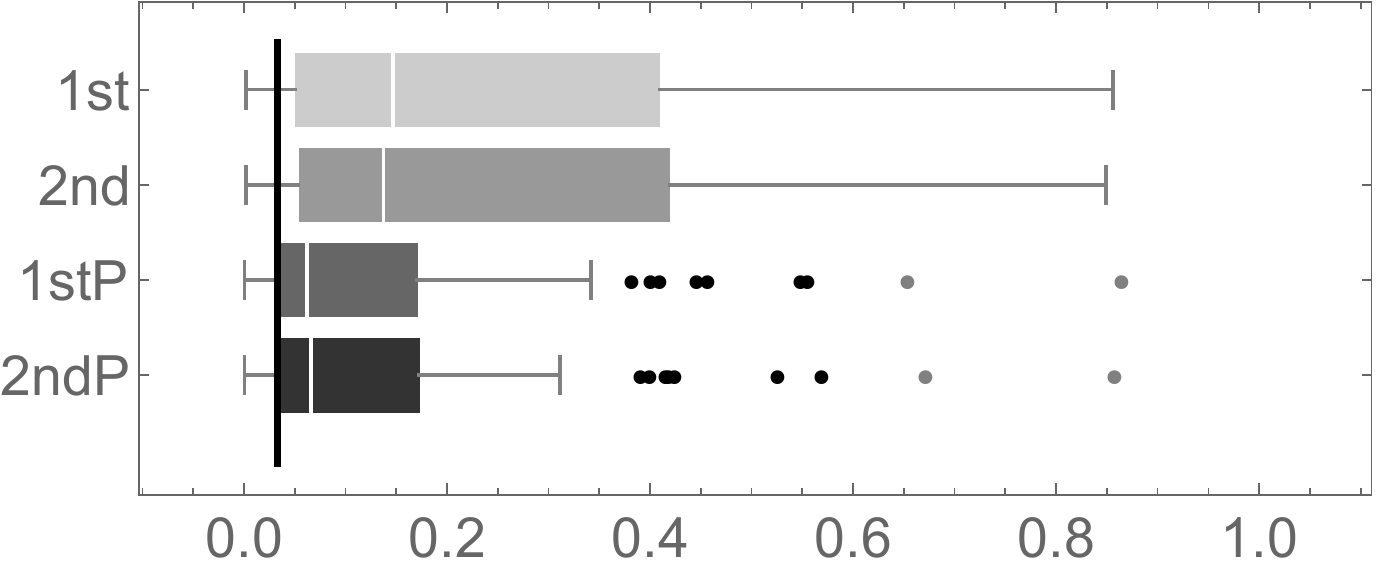}
				\label{fig:13boxM21}
				\vspace{-15pt}
				\caption{$M=2$}
			\end{subfigure}%
			\begin{subfigure}{.45\textwidth}
				\centering
				\includegraphics[width=1\linewidth]{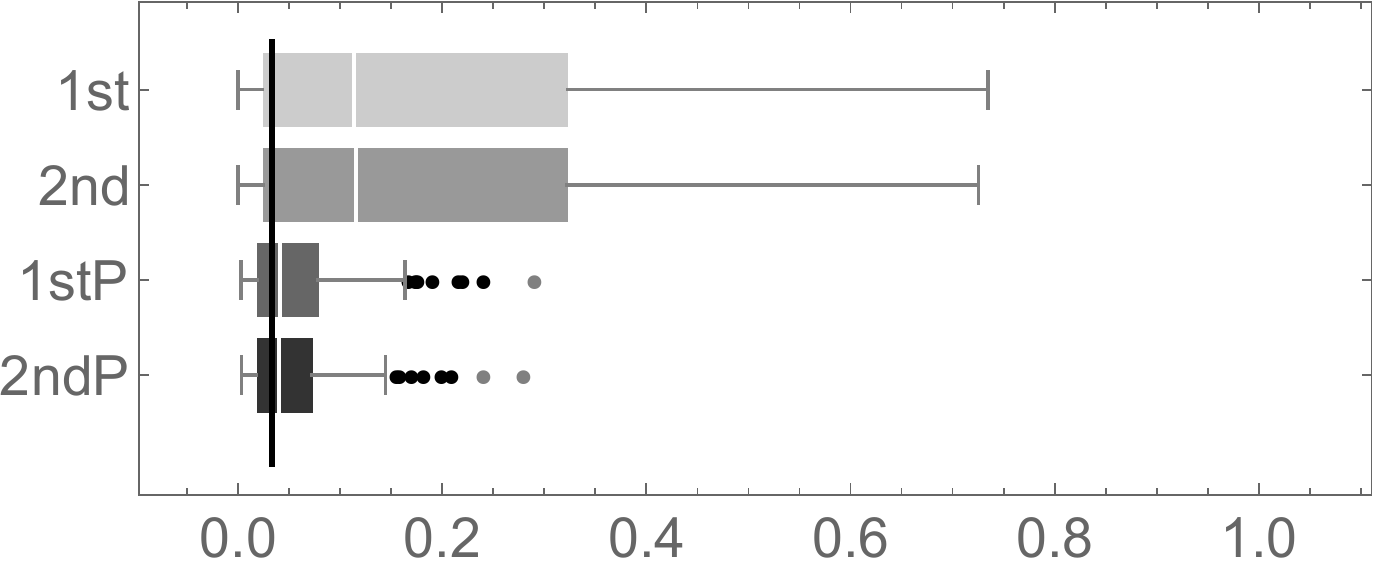}
				\label{fig:13boxM5}
				\vspace{-15pt}
				\caption{$M=5$}
\vspace{10pt}
			\end{subfigure}
			
			\begin{subfigure}{.45\textwidth}
				\centering
				\includegraphics[width=1\linewidth]{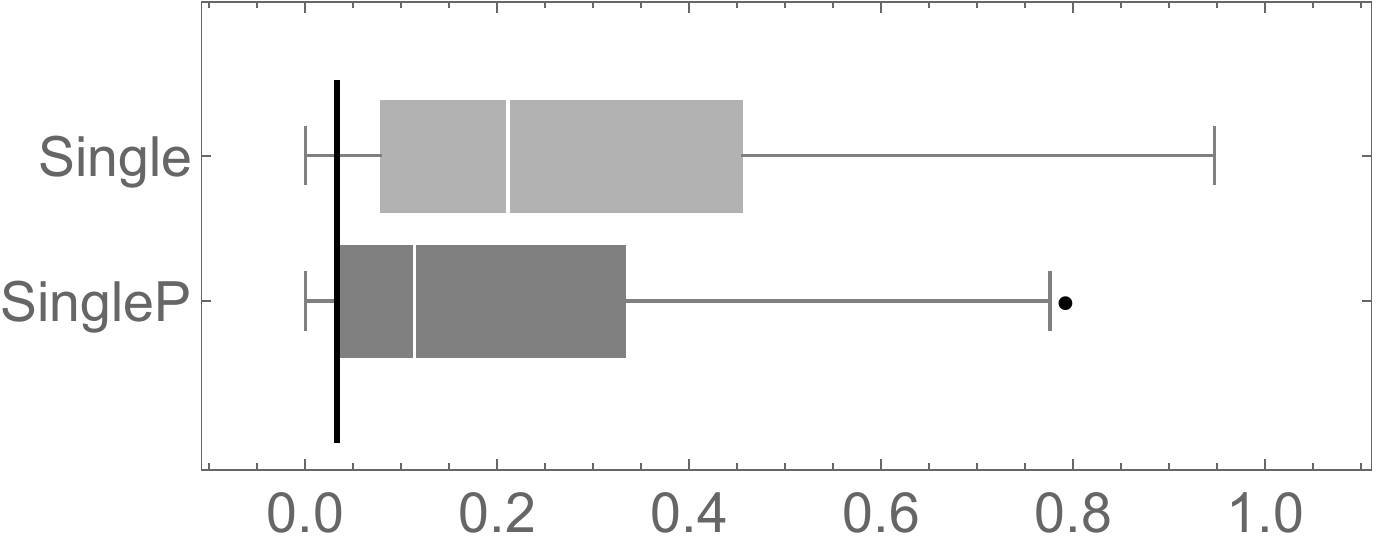}
				\label{fig:13boxSingle}
				\vspace{-15pt}
				\caption{$M=1$}
			\end{subfigure}
			\vspace{-5pt}
			\caption{Box-plots of p-values obtained, when testing the joint significance of A and E, from 100 draws of synthetic datasets using procedures in Section \ref{sec:post} and using Plug-in Sampling method from \cite{moura16}, for $M=1$, $M=2$ and $M=5$ .}\label{fig:box2}
		\end{figure}

		For the second case, we are interested in testing the hypothesis $H_0: \mathbf{AB=C_0}$, where $\mathbf{C_0}$ is a $13\times 13$ matrix consisting of only zeros, with 
		$$\mathbf{A}=\left( \begin{array}{c|c|c}
		\textbf{\Large 0}_{13\times 3} & \textbf{\Large I}_{13} & \textbf{\Large 0}_{13\times 8} \\
		\end{array} \right),$$
		corresponding to the test of joint significance of variables A and E. The p-value obtained for the original data, based on (3) in \cite{moura16}, was $0.033$, thus rejecting their non-significance for $\gamma=0.05$. In Figure \ref{fig:box2}, we can compare the box-plots obtained for the FPPS and Plug-in Sampling methods obtained by generating 100 draws of synthetic datasets, for $M=1$, $M=2$ and $M=5$. The vertical line represents again the original data's p-value.

		From Figure \ref{fig:box2}, we note that the spread of p-values is again larger for our new procedures based on FPPS than the ones from the Plug-in method, majorly leading to a different conclusion from the inference obtained from the original data.

		For the single imputation case, even if the spread of the p-values gathered from the PPS is larger than the ones from the Plug-in Sampling, the distributions of p-values are not that different for the two methods.

		For the two cases studied, the two new FPPS multiple imputation procedures presented have very similar p-values.
		As $M$ increases the spread of the p-values from FPPS becomes smaller and closer to the original data's p-value but at a smaller rate than the p-values from the Plug-in Sampling. 
		
		Nevertheless, this larger spread of the p-values from FPPS will be compensated by an increase of the level of confidentiality, as it can be seen in the next section.
		
		Next, we present the power for the tests
		\begin{equation}\label{eq:Tests}
		\ba{l}\displaystyle H_0:\mathbf{B=B_0(\neq 0)}\; vs\;\; 
		H_1:\mathbf{B=B_1}~~~ \rm{and}\\
		\displaystyle H_0:\mathbf{AB=C_0(\neq 0)}\; vs\;\; H_1:\mathbf{AB=C_1}
		\end{array}
		\end{equation}
		for $\mathbf{B_0}$ equal to $\mathbf{\tilde{\hat{B}}}$, rounded to two decimal places,
		$$\mathbf{A}=\left( \begin{array}{c|c|c}
		\textbf{\Large 0}_{12\times 4} & \textbf{\Large I}_{12} & \textbf{\Large 0}_{12\times 8} \\
		\end{array} \right),$$
		a $12\times 12$ matrix defined appropriately in order to isolate the indicator variables associated with the variable 
		$E$, and $\mathbf{C}_1=\mathbf{AB}_1$ where $\mathbf{B}_1$ takes different values, found in Table 6, with $\mathbf{D}$ a $p\times m$ matrix of $1$'s.

		The power for the synthetic data obtained via FPPS was then simulated as well as the power for the case when these synthetic datasets are treated as if they were the original data. We also simulated the power from the original data and refer to \cite{moura16} for the power values for the synthetic data generated via Plug-in Sampling.

		\begin{table}[h!]
			\footnotesize
			\caption{Power for the tests to the hypothesis (\ref{eq:Tests}), with $\mathbf{B}(1)$, $\mathbf{C}(1)$ and $\mathbf{B}(2)$ and $\mathbf{C}(2)$ denoting the first and second procedures proposed by the authors in Subsections \ref{ssec:mul1} and \ref{ssec:mul2} for FPPS and in \cite{moura16} for Plug-in method.}
			\vspace{-15pt}
			\begin{center}
				\footnotesize\label{ta:power}
				\begin{tabular}{ | l | c | l | c | c | l | c | c | c | l}
					\hline
					Power for       & {orig data} &\multirow{2}{*}{Methods} & {M=1} & \multicolumn{2}{c|}{M=2} & \multicolumn{2}{c|}{M=5} & synt as orig  \\
					$\mathbf{B}_1=$ & $\mathbf{B}$ & & $\mathbf{B}$ & $\mathbf{B}(1)$ & $\mathbf{B}(2)$  & $\mathbf{B}(1)$ & $\mathbf{B}(2)$  & $\mathbf{B}$  \\ \hline \hline
					\multirow{2}{*}{$\mathbf{B}_0+0.005\mathbf{D}$} & \multirow{2}{*}{0.537}  & FPPS & 0.215    & 0.252  & 0.253  & 0.275  & 0.279  & 1.000   \\
					& & Plug-in  & 0.279  & 0.382  & 0.385  & 0.471  & 0.472 & 1.000   \\ \hline
					\multirow{2}{*}{$\mathbf{B}_0*0.95$}   & \multirow{2}{*}{0.945}  & FPPS & 0.535    & 0.634  & 0.637  & 0.700  & 0.700  & 1.000   \\
					&  & Plug-in & 0.679    & 0.840  & 0.841  & 0.906  & 0.909  & 1.000   \\
					\hline
					\hline \hline
					Power for       & {orig data} & \multirow{2}{*}{Methods} & {M=1} & \multicolumn{2}{c|}{M=2} & \multicolumn{2}{c|}{M=5} & synt as orig  \\ 
					$\mathbf{C}_1=$ & $\mathbf{C}$ & & $\mathbf{C}$ & $\mathbf{C}(1)$ & $\mathbf{C}(2)$  & $\mathbf{C}(1)$ & $\mathbf{C}(2)$  & $\mathbf{C}$  \\ \hline \hline
					\multirow{2}{*}{$\mathbf{A(B}_0+3\mathbf{D})$} & \multirow{2}{*}{0.465} &  FPPS & 0.185 & 0.202 & 0.207 & 0.245 & 0.246 & 0.996  \\
					& & Plug-in & 0.284 & 0.334 & 0.343 & 0.416 & 0.418  & 0.975  \\ \hline
					\multirow{2}{*}{$\mathbf{A(B}_0*0.5)$}   & \multirow{2}{*}{0.393} & FPPS & 0.136    & 0.160  & 0.161  & 0.179  & 0.181  & 0.996   \\
					& & Plug-in & 0.197 & 0.271 & 0.279 & 0.326 & 0.327  & 0.959 \\
					\hline	\end{tabular}
				\vspace{-2pt}
			\end{center}
			
		\end{table}
		
		From the power values in Table 6 we may see that tests based on the synthetic data via FPPS show lower values for its power than the ones based in Plug-in generation, as expected, since we are using a method which is supposed to give more confidentiality by generating more perturbed datasets. We may see that these values increase along with the value of $M$, but with a smaller rate than that for Plug-in Sampling, leading to the conclusion that one will need larger values of $M$ to obtain a closer power value to the one registered when testing using the original data. If synthetic data is treated as original, we obtain a larger power than the one obtained for the original data, which is obviously misleading, since the estimated coverage probability will be in fact much smaller than the desired $0.95$.

\section{PRIVACY PROTECTION OF SINGLY VERSUS MULTIPLY IMPUTED SYNTHETIC DATA}
\label{sec:pri}
		
		In order to evaluate the level of protection and at the same time compare it with the level obtained from synthetic data generated via Plug-in Sampling, we perform, in this section, a similar evaluation as in \cite{moura16} using CPS data. Let us consider $\mathbf{W}_l=(\mathbf{w}_{1l},...,\mathbf{w}_{nl})$, $l=1,...,M$, $M$ synthetic datasets generated via FPPS, where $\mathbf{w}_{il}=(w_{1il},...,w_{mil})', i=1,...,n$. The estimate of the original values $\mathbf{y}_i=(y_{1i},...,y_{mi})'$ will be $\hat{\mathbf{y}}_i=\frac{1}{M}\sum_{l=1}^M \mathbf{w}_{il}$. Let us recall the three criteria used in \cite{moura16} as measures of the level of privacy protection:
		\begin{equation}\label{eq:conf}
		\begin{array}{c}\displaystyle
		\Gamma_{1,\epsilon}=\frac{1}{m n}\sum_{j=1}^m\sum_{i=1}^n 
		Pr\left[\,\left|\frac{\hat{y_{ji}}-y_{ji}}{y_{ji}}\right|<\epsilon\,\Bigl|\,\mathbf{Y}\right];\medskip\\
		\displaystyle
		\Gamma_{2,\epsilon}=\frac{1}{n}\sum_{i=1}^n 
		Pr\left[\sqrt{\frac{1}{m}\sum_{j=1}^m\frac{(\hat{y_{ji}}-y_{ji})^2}{y_{ji}^2}}<\epsilon\,\Bigl|\,\mathbf{Y}\right];\medskip\\
		\displaystyle
		\Gamma_{3,\epsilon}=Pr\left[\frac{1}{m 
			n}\sum_{j=1}^m\sum_{i=1}^n\left|\frac{\hat{y_{ji}}-y_{ji}}{y_{ji}}\right|<\epsilon\,\Bigl|\,\mathbf{Y}\right].
		\end{array}
		\end{equation}
		
		Let us also consider, from $\Gamma_{1,\epsilon}$, the following quantity, for $i=1,...n$ and $j=1,..,m$,
		$$D_{1,\epsilon,ji}=Pr\left[\left|\frac{\hat{y_{ji}}-y_{ji}}{y_{ji}}\right|<\epsilon\,\Bigl|\,\mathbf{Y}\right]$$
		and, from $\Gamma_{3,\epsilon}$,
		$$D_3=\frac{1}{m n}\sum_{j=1}^m\sum_{i=1}^n\left|\frac{\hat{y_{ji}}-y_{ji}}{y_{ji}}\right|.$$ We use a Monte Carlo simulation with 
		$10^4$ iterations to estimate all three measures in (\ref{eq:conf}) based on the $n=141$ households in the CPS data. In Table 7, we show the 
		values of $\Gamma_{1,0.01}$, $\Gamma_{2,0.01}$ and the minimum, 1st quartile ($Q_1$), median, 3rd quartile ($Q_3$) and maximum of 
		$D_{1,\epsilon}$, displaying also the values gathered when using Plug-in Sampling. In Table 8, we show the values of $\Gamma_{3,0.1}$ and the minimum, $Q_1$, median, $Q_3$ and maximum of $D_3$ also displaying the values gathered when using Plug-in Sampling. 
		
		\begin{table}[h!]
			\caption{Values of $\Gamma_{1,0.01}$, $\Gamma_{2,0.01}$ and a summary of the distribution of $D_{1,0.01}$.}
			\footnotesize \label{table:t12}
			\vspace{-20pt}
			\begin{center}
				\begin{tabular}{ | c | l | c | c | c | c | c | c | c | }\hline 
					$M$ & Method &	$\Gamma_{1,0.01}$ & $\Gamma_{2,0.01}$ & Min & $Q_1$ & Median & $Q_3$ & Max  \\ \hline \hline
					\multirow{2}{*}{$M=1$} & FPPS & 0.0602  & 0.0005 & 0 & 0.0385 & 0.0507 & 0.0784 & 0.1455 \\ 
					& Plug-in & 0.0631  & 0.0006 & 0 & 0.0398 & 0.0552 & 0.0854 & 0.1491 \\ \hline
					\multirow{2}{*}{$M=2$} & FPPS & 0.0702  & 0.0009 & 0 & 0.0357 & 0.0624 & 0.0910 & 0.1945 \\ 
					&Plug-in & 0.0754  & 0.0010 & 0 & 0.0331 & 0.0697 & 0.0954 & 0.2134 \\ \hline
					\multirow{2}{*}{$M=5$} & FPPS & 0.0797  & 0.0012 & 0 & 0.0214 & 0.0711 & 0.1136 & 0.2785 \\ 
					&Plug-in & 0.0879  & 0.0018 & 0 & 0.0110 & 0.0792 & 0.1284 & 0.3279 \\ \hline
				\end{tabular}
			\end{center}
		\end{table}
		
		\begin{table}[h!] 
			\caption{Values of $\Gamma_{3,0.1}$ and a summary of the distribution of $D_3$.}
			\footnotesize\label{table:t3}
			\vspace{-20pt}
			\begin{center}
				\begin{tabular}{| c | l | c | c | c | c | c | c |  }\hline 
					$M$    & Method &$\Gamma_{3,0.1}$& Min & $Q_1$ & Median & $Q_3$ & Max  \\ \hline \hline
					\multirow{2}{*}{$M=1$} & FPPS & 0.0000 & 0.1091  & 0.1248 & 0.1287 & 0.1325 & 0.1544 \\ [0pt]
					&Plug-in  &0.0000 &  0.1050  & 0.1202 & 0.1233 & 0.1264 & 0.1379 \\ [0pt] \hline
					\multirow{2}{*}{$M=2$} & FPPS & 0.0021 & 0.0960  & 0.1088 & 0.1116 & 0.1145 & 0.1324 \\ [0pt]
					&Plug-in & 0.0694 & 0.0948  & 0.1026 & 0.1051 & 0.1072 & 0.1159 \\ [0pt] \hline
					\multirow{2}{*}{$M=5$} & FPPS &0.5008      & 0.0896  & 0.0980 & 0.1000 & 0.1020 & 0.1131 \\
					&Plug-in & 1.0000 & 0.0846  & 0.0905 & 0.0920 & 0.0936 & 0.0992 \\ \hline
				\end{tabular}
				\vspace{-10pt}
			\end{center}
		\end{table}

		Looking at Tables \ref{table:t12} and \ref{table:t3}, we observe that the values of the privacy measures in (\ref{eq:conf}) increase for increasing values of $M$ for both procedures developed in Subsections \ref{ssec:mul1} and \ref{ssec:mul2}, showing that the disclosure risk increases with the increase in the number of released synthetic datasets. Compared with the measures obtained under Plug-in Sampling, we may observe a smaller disclosure risk in all cases, leading to the conclusion that the proposed FPPS procedures have an overall higher level of confidentiality. Regarding measures $\Gamma_{2,\epsilon}$ and $\Gamma_{3,\epsilon}$ this increase reaches in some cases an increase of 50\% or more in confidentiality. In the single imputation case, under the PPS we also register an increase of confidentiality when comparing the same measure under Plug-in Sampling, nevertheless this increase is relatively small.

		\section{CONCLUDING REMARKS}
		\label{sec:con}
		
		In this paper the authors derive likelihood-based exact inference for single and multiple imputation cases where synthetic datasets are generated via Fixed-Posterior Predictive Sampling (FPPS).  If only one synthetic dataset is released, then FPPS is equivalent to the usual Posterior Predictive Sampling (PPS) method.  Thus the proposed methodology can be used to analyze a singly imputed synthetic data set generated via PPS under the multivariate linear regression (MLR) model.  Therefore this work fills a gap in the literature because the state of the art methods apply only to multiply imputed synthetic data.  Under the MLR model, the authors derived two different exact inference procedures for the matrix of regression coefficients, when multiply imputed synthetic datasets are released. It is shown that the methodologies proposed lead to confidence sets matching the expected level of confidence, for all sample sizes. 
		Furthermore, while the second proposed procedure displays a better precision for smaller samples and/or smaller values of $M$ by yielding smaller confidence sets, the two procedures concur for larger sample sizes and larger values of $M$, as it is corroborated in theory by remarks 2.2 and 2.3.
		When compared with inference procedures for Plug-in Sampling, the procedures proposed based on FPPS lead to synthetic datasets that give respondents a higher level of confidentiality, that is, a reduced disclosure risk, nevertheless at the expense of accuracy, since the confidence sets are larger, as illustrated in the application with the CPS data. Once likelihood-based exact inferential methods are now made available both for FPPS/PPS and Plug-in Sampling, it is therefore the responsibility of those in charge of releasing the data to decide which method to use in order to better respect the demands and objectives of their institution.
		
\begin{acknowledgments}
		
		Ricardo Moura's research is supported by a Fulbright Research Grant, and he sincerely thanks the faculty of Mathematics and Statistics at UMBC for their support and encouragement. Ricardo Moura and Carlos A. Coelho also thank FCT (Portuguese Foundation for Science and Technology) project UID/MAT/00297/2013 awarded through CMA/UNL. Martin Klein and Bimal Sinha thank Laura McKenna, Eric Slud, William Winkler, and 
		Tommy Wright at the U.S. Census Bureau for their support.
		The authors would also like to thank the referees for the helpful comments and suggestions leading to the improvement of the paper.
\end{acknowledgments}

		\renewcommand{\theequation}{\thesection.\arabic{equation}}
		\setcounter{equation}{0}
		
		\begin{appendices}
			
			\section{Proof of Theorems 2.1 and 2.2 and Corollaries \ref{cor:pdf} and \ref{cor:dist}} \label{App:A}

\begin{proof}[Proof of Theorem \ref{thm:pdf}] 
Given $(\mathbf{\tilde{B}},\mathbf{\tilde{\Sigma}})$, from (\ref{eq:synt}) we have that, for every $j=1,...,M$, 
				$$\mathbf{W}_j'|_{\mathbf{\tilde{B},\tilde{\Sigma}}}\sim N_{nm}(\mathbf{X'\tilde{B}},\mathbf{\tilde{\Sigma}\otimes I}_n)\implies 
				\mathbf{B}^{\bullet}_j|_{\mathbf{\tilde{B},\tilde{\Sigma}}}\sim N_{pm}(\mathbf{\tilde{B}},\mathbf{\tilde{\Sigma}} 
				\otimes (\mathbf{XX'})^{-1})$$
				and
				$$(n-p)\mathbf{S}_j^{\bullet}|_{\mathbf{\tilde{\Sigma}}}\sim W_m(\mathbf{\tilde{\Sigma}},n-p).$$
				Therefore, we have for $\mathbf{\overline{B}}_M^\bullet$ and $\mathbf{\overline{S}}_M^\bullet$ in (\ref{eq:parameters1st}),
				$$\mathbf{\overline{B}}_M^\bullet|_{\mathbf{\tilde{B},\tilde{\Sigma}}}=\frac{1}{M}\sum_{j=1}^M \mathbf{B}^\bullet_j|_{\mathbf{\tilde{B},\tilde{\Sigma}}}\sim N_{pm}\left(\tilde{\mathbf{B}},\frac{1}{M}\tilde{\mathbf{\Sigma}}\otimes (\mathbf{XX'})^{-1}\right)$$
				and
				$$M(n-p)\mathbf{\overline{S}}^\bullet_M|_{\mathbf{\tilde{\Sigma}}}=(n-p)\sum_{j=1}^M\mathbf{\overline{S}}^\bullet_j|_{\mathbf{\tilde{\Sigma}}}\sim W_m(\mathbf{\tilde{\Sigma}},M(n-p)).$$
				
				Since $\mathbf{\overline{B}}_M^\bullet$ and $\mathbf{\overline{S}}_M^\bullet$ are independent, the conditional joint pdf of $(\mathbf{\overline{B}}^{\bullet}_M,\mathbf{\overline{S}}^{\bullet}_M)$, given $\tilde{\mathbf{B}}$ and $\tilde{\mathbf{\Sigma}}$, is
\begin{equation}
\label{eq:syntjoint}
\ba{l}
\small f(\mathbf{\overline{B}}^{\bullet}_M,\mathbf{\overline{S}}^{\bullet}_M\vert \mathbf{\tilde{B}},\mathbf{\tilde{\Sigma}}) \propto \\
\hskip 1.55cm 
e^{-\frac{1}{2}tr\lbrace M\mathbf{\tilde{\Sigma}}^{-1}\left[(\mathbf{\overline{B}}^{\bullet}_M-\mathbf{\tilde{B}})'\mathbf{XX'}(\mathbf{\overline{B}}^{\bullet}_M-\mathbf{\tilde{B}})+M(n-p)\mathbf{\overline{S}}^{\bullet}_M\right]\rbrace}\times 
				\frac{\vert\mathbf{\overline{S}}^{\bullet}_M\vert^{\frac{M(n-p)-m-1}{2}}}{\vert\mathbf{\tilde{\Sigma}}\vert^{\frac{M(n-p)+p}{2}}}, 
\ea
\end{equation}
while, due to the independence of $ \mathbf{\tilde{\Sigma}^{-1}} $ and $ \mathbf{\tilde{B}} $, generated from (\ref{eq:postSigma}) and (\ref{eq:postB}), respectively, the joint pdf of $(\mathbf{\tilde{B}},\mathbf{\tilde{\Sigma}^{-1}})$, given $\mathbf{S}$, is
\begin{equation}\label{eq:postdist}
f(\mathbf{\tilde{B}},\mathbf{\tilde{\Sigma}^{-1}}|\mathbf{S})\propto |\mathbf{\tilde{\Sigma}}|^{-p/2} e^{-\frac{1}{2}tr\lbrace\mathbf{\tilde{\Sigma}}^{-1}\left[(\mathbf{\tilde{B}-\hat{B}})'\mathbf{XX'}(\mathbf{\tilde{B}-\hat{B}})+(n-p)\mathbf{S}\right]\rbrace}\frac{\vert \mathbf{S}\vert^{\frac{n+\alpha-p-m-1}{2}}}{\vert\mathbf{\tilde{\Sigma}}\vert^{\frac{n+\alpha-p}{2}-m-1}}.\!\!\!
				\end{equation}
				On the other hand, given the independence of $\mathbf{\hat{B}}$ and $\mathbf{S}$, defined in (\ref{eq:orB}) and (\ref{eq:orS}), the joint pdf of $(\mathbf{\hat{B}},\mathbf{S})$ is given by
				\begin{equation}\label{eq:origjoint} f(\mathbf{\hat{B}},\mathbf{S})\propto
				e^{-\frac{1}{2}tr\lbrace\mathbf{\Sigma}^{-1}\left[(\mathbf{\hat{B}-B})'\mathbf{XX'}(\mathbf{\hat{B}-B})+(n-p)\mathbf{S}\right]\rbrace}\frac{\vert 
					\mathbf{S}\vert^{\frac{n-p-m-1}{2}}}{\vert\mathbf{\Sigma}\vert^{\frac{n}{2}}}. \end{equation}

				Thus, by multiplying the three pdf's in (\ref{eq:syntjoint}), (\ref{eq:postdist}) and (\ref{eq:origjoint}), we obtain the joint pdf of 
				$(\mathbf{\overline{B}}^{\bullet}_M,\mathbf{\overline{S}}^{\bullet}_M,\mathbf{\tilde{B},\tilde{\Sigma}^{-1},\hat{B},S})$.
				
				Since
				$$tr\lbrace M(\mathbf{\overline{B}}^{\bullet}_M-\mathbf{\tilde{B}})'\mathbf{XX'}(\mathbf{\overline{B}}^{\bullet}_M-\mathbf{\tilde{B}})\rbrace=tr\lbrace M(\mathbf{\tilde{B}}-\mathbf{\overline{B}}^{\bullet}_M)'\mathbf{XX'}(\mathbf{\tilde{B}}-\mathbf{\overline{B}}^{\bullet}_M)\rbrace,$$
				and since from Appendix \ref{Aapp:last2} we may write
				$$\hspace{-100pt}M(\mathbf{\tilde{B}}-\mathbf{\overline{B}}^{\bullet}_M)'\mathbf{XX'}(\mathbf{\tilde{B}}-\mathbf{\overline{B}}^{\bullet}_M)+(\mathbf{\tilde{B}-\hat{B}})'\mathbf{XX'}(\mathbf{\tilde{B}-\hat{B}})=$$
				$$\hspace{-50pt}=\!(M+1)\!\left[\mathbf{\tilde{B}}\!-\!\frac{1}{M+1}\mathbf{(B^{\bullet}+\hat{B})}\right]'\!\mathbf{XX'}\left[\mathbf{\tilde{B}}\!-\!\frac{1}{M+1}\mathbf{(B^{\bullet}+\hat{B})}\right]\!$$
				$$\hspace{220pt}+\frac{M}{M+1}(\mathbf{B^{\bullet}\!-\!\hat{B}})'\mathbf{XX}'(\mathbf{B^{\bullet}\!-\!\hat{B}}),$$
by integrating out $\mathbf{\tilde{B}}$, we obtain the joint pdf of $(\mathbf{\overline{B}}^{\bullet}_M,\mathbf{\overline{S}}^{\bullet}_M,\mathbf{\tilde{\Sigma}^{-1},\hat{B},S})$ proportional to
				\begin{flalign} \nonumber
				e^{-\frac{1}{2}tr\lbrace\mathbf{\tilde{\Sigma}}^{-1}\left[\frac{M}{M+1}(\mathbf{\overline{B}}^{\bullet}_M-\mathbf{\hat{B}})'\mathbf{XX'}(\mathbf{\overline{B}}^{\bullet}_M-\mathbf{\hat{B}})+(n-p)(M\mathbf{\overline{S}}^{\bullet}_M+\mathbf{S})\right]+\mathbf{\Sigma}^{-1}\left[(\mathbf{\hat{B}-B})'\mathbf{XX'}(\mathbf{\hat{B}-B})+(n-p)\mathbf{S}\right]\rbrace}\\ 
\nonumber				\times \frac{\vert\mathbf{\overline{S}}^{\bullet}_M\vert^{\frac{M(n-p)-m-1}{2}}}{\vert\mathbf{\tilde{\Sigma}}\vert^{\frac{M(n-p)+n-\alpha}{2}-m-1}}\frac{\vert \mathbf{S}\vert^{n+\frac{\alpha}{2}-p-m-1}}{\vert\mathbf{\Sigma}\vert^{\frac{n}{2}}}.
				\end{flalign}
				
				Since 
				\begin{align*} 
				tr\left\{\frac{M}{M+1}\mathbf{\tilde{\Sigma}^{-1}}(\mathbf{\overline{B}}^{\bullet}_M-\mathbf{\hat{B}})'(\mathbf{XX'})(\mathbf{\overline{B}}^{\bullet}_M-\mathbf{\hat{B}})+\mathbf{\Sigma^{-1}}(\mathbf{\hat{B}-B})'(\mathbf{XX'})(\mathbf{\hat{B}-B})\right\}=\\ 
				tr\left\{\mathbf{XX'}\left[\frac{M}{M+1}(\mathbf{\overline{B}}^{\bullet}_M-\mathbf{\hat{B}})\mathbf{\tilde{\Sigma}^{-1}}(\mathbf{\overline{B}}^{\bullet}_M-\mathbf{\hat{B}})'+(\mathbf{\hat{B}-B})\mathbf{\Sigma^{-1}}(\mathbf{\hat{B}-B})'\right]\right\}
				\end{align*} 
				and since from the identities in 1.-3. in Appendix B1 in \cite{moura16} we may write
				\begin{align*} & 
				\frac{M}{M+1}(\mathbf{\overline{B}}^{\bullet}_M-\mathbf{\hat{B}})\mathbf{\tilde{\Sigma}^{-1}}(\mathbf{\overline{B}}^{\bullet}_M-\mathbf{\hat{B}})'+(\mathbf{\hat{B}-B})\mathbf{\Sigma^{-1}}(\mathbf{\hat{B}-B})'=\\ 
				&=\left[\mathbf{\hat{B}}-\left(\frac{M}{M+1}\mathbf{\overline{B}}^{\bullet}_M\mathbf{\tilde{\Sigma}^{-1}}+\mathbf{B}\mathbf{\Sigma^{-1}}\right)\left(\frac{M}{M+1}\mathbf{\tilde{\Sigma}^{-1}}+\mathbf{\Sigma^{-1}}\right)^{-1}\right]\\ 
				&\left(\!\frac{M}{M+1}\!\mathbf{\tilde{\Sigma}^{-1}}\!+\!\mathbf{\Sigma^{-1}}\right)\left[\mathbf{\hat{B}}\!-\!\left(\frac{M}{M+1}\mathbf{\overline{B}}^{\bullet}_M\mathbf{\tilde{\Sigma}^{-1}}\!+\!\mathbf{B}\mathbf{\Sigma^{-1}}\right)\left(\frac{M}{M+1}\mathbf{\tilde{\Sigma}^{-1}}\!+\!\mathbf{\Sigma^{-1}}\right)^{-1}\right]'\\ 
				&\;\;\;\;\;\;\;\;\;\;\;\;\;\;\;\;\;\;\;\;\;\;\;\;\;\;\;\;\;\;\;\;\;\;\;\;\;\;\;\;\;\;\;\;\;\;\;\;\;\;\;\;\;+(\mathbf{\overline{B}}^{\bullet}_M-\mathbf{B})\left(\frac{M+1}{M}\mathbf{\tilde{\Sigma}}+\mathbf{\Sigma}\right)^{-1}(\mathbf{\overline{B}}^{\bullet}_M-\mathbf{B})', \end{align*}
				integrating out $\mathbf{\hat{B}}$ we will have the joint pdf of $(\mathbf{\overline{B}}^{\bullet}_M,\mathbf{\overline{S}}^{\bullet}_M,\mathbf{\tilde{\Sigma}^{-1},S})$ proportional to
				\begin{flalign*} \nonumber
				& e^{-\frac{1}{2}tr\lbrace(\frac{M+1}{M}\mathbf{\tilde{\Sigma}+\Sigma)}^{-1}(\mathbf{\overline{B}}^{\bullet}_M-\mathbf{\hat{B}})'\mathbf{XX'}(\mathbf{\overline{B}}^{\bullet}_M-\mathbf{\hat{B}})+(n-p)\mathbf{\tilde{\Sigma}^{-1}}(M\mathbf{\overline{S}}^{\bullet}_M+\mathbf{S})+(n-p)\mathbf{\Sigma}^{-1}\mathbf{S}\rbrace}\nonumber\\ 
				&~~~~~~~~~~~~~~~~~~~~~~\times \frac{\vert\mathbf{\overline{S}}^{\bullet}_M\vert^{\frac{M(n-p)-m-1}{2}}}{\vert\mathbf{\tilde{\Sigma}}\vert^{\frac{M(n-p)+n-\alpha}{2}-m-1}}\frac{\vert \mathbf{S}\vert^{n+\frac{\alpha}{2}-p-m-1}}{\vert\mathbf{\Sigma}\vert^{\frac{n}{2}}} \left|\frac{M}{M+1}\mathbf{\tilde{\Sigma}^{-1}}+\mathbf{\Sigma^{-1}}\right|^{-p/2}.
				\end{flalign*} Consequently, if we integrate out $\mathbf{S}$ we will end up with the joint pdf of \linebreak $(\mathbf{\overline{B}}^{\bullet}_M,\mathbf{\overline{S}}^{\bullet}_M,\mathbf{\tilde{\Sigma}^{-1}})$ proportional to
				\begin{flalign}\label{eq:likelihood}
				& e^{-\frac{1}{2}tr\lbrace(\frac{M+1}{M}\mathbf{\tilde{\Sigma}+\Sigma)}^{-1}(\mathbf{\overline{B}}^{\bullet}_M-\mathbf{B})'\mathbf{XX'}(\mathbf{\overline{B}}^{\bullet}_M-\mathbf{B})+M(n-p)\mathbf{\tilde{\Sigma}^{-1}}\mathbf{\overline{S}}^{\bullet}_M\rbrace}\\ 
				&{\small \times \frac{\vert\mathbf{\overline{B}}^{\bullet}_M\vert^{\frac{M(n-p)-m-1}{2}}}{\vert\mathbf{\tilde{\Sigma}}\vert^{\frac{M(n-p)+n-\alpha}{2}-m-1}}\vert\mathbf{\Sigma}\vert^{-\frac{n}{2}} \left|\frac{M}{M+1}\mathbf{\tilde{\Sigma}^{-1}}+\mathbf{\Sigma^{-1}}\right|^{-p/2} |\mathbf{\tilde{\Sigma}^{-1}+\Sigma^{-1}}|^{-\frac{2n+\alpha-2p-m-1}{2}}}\nonumber
				\end{flalign}
				as we wanted to prove. It is easy to see that in (\ref{eq:likelihood}), $\mathbf{\overline{S}}^{\bullet}_M$ and 
				$\mathbf{\overline{B}}^{\bullet}_M$, given $\mathbf{\tilde{\Sigma}}^{-1}$, are separable, with the distributions in the body of the Theorem.\end{proof}

			\begin{proof}[Proof of Theorem \ref{thm:dist}]  
				From the distributions of $\mathbf{\overline{S}}^{\bullet}_M$ and $\mathbf{\overline{B}}^{\bullet}_M$ in Theorem \ref{thm:pdf}, 
				and by Theorem 2.4.1 in \cite{kollo05} we have 
				that, for $p\geq m,$
				$$(\mathbf{\overline{B}}^{\bullet}_M-\mathbf{B})'(XX')(\mathbf{\overline{B}}^{\bullet}_M-\mathbf{B})|\mathbf{_{\tilde{\Sigma}^{-1}}}\sim W_m\left(\frac{M+1}{M}\mathbf{\tilde{\Sigma}+\Sigma},p\right).$$
				From Theorem 2.4.2 in \cite{kollo05} and Subsection 7.3.3 in \cite{anderson84} we have
				\begin{equation}\label{eq:H}\small
				\mathbf{H}=\left(\frac{M+1}{M}\mathbf{\tilde{\Sigma}\!+\!\Sigma}\right)^{\!\!-\frac{1}{2}}\!(\mathbf{\overline{B}}^{\bullet}_M\!-\!\mathbf{B)'(XX')}(\mathbf{\overline{B}}^{\bullet}_M\!-\!\mathbf{B)}\left(\frac{M+1}{M}\mathbf{\tilde{\Sigma}\!+\!\Sigma}\right)^{\!\!\prime -\frac{1}{2}}
\sim W_m(\mathbf{I},p)
				\end{equation}
				and
				\begin{equation}\label{eq:G}\mathbf{G}=M(n-p)\mathbf{\tilde{\Sigma}}^{-\frac{1}{2}}\mathbf{\overline{S}}^{\bullet}_M\mathbf{\tilde{\Sigma}}'^{-\frac{1}{2}}
				\sim W_m(\mathbf{I},M(n-p)).
				\end{equation}
				We may thus write $T^\bullet_M$ in (\ref{eq:T1st}) as 
				$$T^\bullet_M=\frac{|(\mathbf{\overline{B}}^{\bullet}_M-\mathbf{B})'(XX')(\mathbf{\overline{B}}^{\bullet}_M-\mathbf{B})|}{|M(n-p)\mathbf{\overline{S}}^{\bullet}_M|}=\frac{\left|\frac{M+1}{M}\mathbf{\tilde{\Sigma}+\Sigma}\right|}{|\mathbf{\tilde{\Sigma}}|}\times\frac{|\mathbf{H}|}{|\mathbf{G}|},$$ 
				where, 
$|\mathbf{G}|\sim 
				\prod_{i=1}^{m}\chi^2_{n-p-i+1}$ and $|\mathbf{H}| \sim \prod_{i=1}^{m}\chi^2_{p-i+1}$, with independent chi-square random variables in each product, we end up with a product of independent F-distributions, due to the independence of $\mathbf{H}$ and $\mathbf{G}$, inherited from the independence of $\mathbf{\overline{B}}^{\bullet}_M$ and $\mathbf{\overline{S}}^{\bullet}_M$. So, conditionally on $\mathbf{\tilde{\Sigma}^{-1}}$, we have
				$$T^\bullet_M|_\mathbf{\tilde{\Sigma}^{-1}}\sim \left\{\prod_{i=1}^{m}\frac{p-i+1}{M(n-p)-i+1} F_{p-i+1,n-p-i+1}\right\}\times \left|\mathbf{\tilde{\Sigma}}^{-1}\left( \frac{M+1}{M}\mathbf{\tilde{\Sigma}+\Sigma }\right)\right|,$$
where
$$\left|\mathbf{\tilde{\Sigma}^{-1}}\left(\frac{M+1}{M}\mathbf{\tilde{\Sigma}}+\mathbf{\Sigma}\right)\right|=\left|\frac{M+1}{M}\mathbf{I+\mathbf{\tilde{\Sigma}^{-1}\Sigma}}\right|=\left|\frac{M+1}{M}\mathbf{\Sigma^{-1}+\mathbf{\tilde{\Sigma}^{-1}}}\right||\mathbf{\Sigma}|$$
				$$~~~~~~~~~=\left|\mathbf{\Sigma}^{1/2}\right|\left|\frac{M+1}{M}\mathbf{\Sigma^{-1}+\mathbf{\tilde{\Sigma}^{-1}}}\right|\left|\mathbf{\Sigma}^{1/2}\right|=\left|\frac{M+1}{M}\mathbf{I+\Sigma^{1/2}\tilde{\Sigma}^{-1}\Sigma^{1/2}}\right|.$$
				
				As such, from (\ref{eq:likelihood}), integrating out $\mathbf{\overline{B}}^{\bullet}_M$ and $\mathbf{\overline{S}}^{\bullet}_M$, we end up with the pdf of $\mathbf{\tilde{\Sigma}}^{-1}$ proportional to
				\begin{flalign*}
&|\mathbf{\tilde{\Sigma}}|^{\frac{M(n-p)}{2}}\left|\frac{M+1}{M}\mathbf{\tilde{\Sigma}+\Sigma}\right|^{\frac{p}{2}}\frac{1}{\vert\mathbf{\tilde{\Sigma}}\vert^{\frac{M(n-p)+n-\alpha}{2}-m-1}}\vert\mathbf{\Sigma}\vert^{-\frac{n}{2}}\medskip\\
				& ~~~~~~~~~~~~~~~~~~~~~~~~~~~~~~~~~\times \left|\frac{M}{M+1}\mathbf{\tilde{\Sigma}^{-1}}+\mathbf{\Sigma^{-1}}\right|^{-p/2} |\mathbf{\tilde{\Sigma}^{-1}+\Sigma^{-1}}|^{-\frac{2n+\alpha-2p-m-1}{2}}\medskip\\[10pt]
&=|\mathbf{\tilde{\Sigma}}^{-1}|^{\frac{n+\alpha-2m-2}{2}}\left|\frac{M+1}{M}\mathbf{\tilde{\Sigma}+\Sigma}\right|^{\frac{p}{2}}\vert\mathbf{\Sigma}\vert^{-\frac{n}{2}}\\
				& ~~~~~~~~~~~~~~~~~~~~~~~~~~~~~~~~~\times \left|\frac{M}{M+1}\mathbf{\tilde{\Sigma}^{-1}}+\mathbf{\Sigma^{-1}}\right|^{-p/2} |\mathbf{\tilde{\Sigma}^{-1}+\Sigma^{-1}}|^{-\frac{2n+\alpha-2p-m-1}{2}}.
				\end{flalign*}

				Making the transformation $\mathbf{\Omega}\!=\!\mathbf{\Sigma^{\frac{1}{2}}\tilde{\Sigma}^{-1}\Sigma^{\frac{1}{2}}}$, which implies $\mathbf{\tilde{\Sigma}^{-1}}\!=\!\mathbf{\Sigma^{-\frac{1}{2}} \Omega \Sigma^{-\frac{1}{2}}}$, with the Jacobian of the transformation from $\mathbf{\tilde{\Sigma}}^{-1}$ to $\mathbf{\Omega}$ being $|\mathbf{\Sigma}|^{-\frac{m+1}{2}}$, we have the pdf of $\mathbf{\Omega}$ proportional to
				\begin{flalign*}
				|\mathbf{\Omega}|^{\frac{n+\alpha-2m-2}{2}}\left|\frac{M+1}{M}\mathbf{\Omega^{-1}+I_m}\right|^{\frac{p}{2}}\left|\frac{M}{M+1}\mathbf{\Omega}+\mathbf{I_m}\right|^{-p/2} |\mathbf{\Omega+I_m}|^{-\frac{2n+\alpha-2p-m-1}{2}}.
				\end{flalign*}
				Since $|\frac{M+1}{M}\mathbf{\Omega^{-1}+I_m}|^{\frac{p}{2}}=\left(\frac{M+1}{M}\right)^{p/2}|\frac{M}{M+1}\mathbf{\Omega+I_m}|^{\frac{p}{2}}|\mathbf{\Omega}|^{-\frac{p}{2}}$ we end up with
				\begin{flalign*}
				f(\mathbf{\Omega})\propto|\mathbf{\Omega}|^{\frac{n+\alpha-p-2m-2}{2}}\times |\mathbf{\Omega+I_m}|^{-\frac{2n+\alpha-2p-m-1}{2}}
				\end{flalign*}
				independent of $\mathbf{\Sigma}$. Therefore, we may conclude that
				$$T^\bullet_M|_\mathbf{\Omega}\sim \left\{\prod_{i=1}^{m}\frac{p-i+1}{n-p-i+1} F_{p-i+1,M(n-p)-i+1}\right\} \left|\frac{M+1}{M}\mathbf{I_m+\Omega}\right|$$
				where from \cite[Theorem 8.2.8.]{muirhead05} $\mathbf{\Omega}$ has the same distribution as $\mathbf{A}_1^{\frac{1}{2}}\mathbf{A}_2^{-1}\mathbf{A}_1^{\frac{1}{2}}$ with $\mathbf{A}_1\sim W_m(\mathbf{I}_m, n+\alpha-p-m-1)$ and $\mathbf{A}_2\sim W_m(\mathbf{I}_m, n-p)$, two independent random variables.
			\end{proof}

			\begin{proof}[Proof of Corollary \ref{cor:pdf}]
				
				The proof is identical to the proof of Theorem \ref{thm:pdf} replacing the 
				joint pdf of $(\mathbf{\overline{B}}^{\bullet}_M,\mathbf{\overline{S}}^{\bullet}_M)$ by the joint pdf of $(\mathbf{\overline{B}}^{\bullet}_M,\mathbf{S}^\bullet_{comb})$, noting that we 
				have
				\[(Mn-p)\mathbf{S}^\bullet_{comb}|_{\tilde{\mathbf{\Sigma}}}\sim W_m(\mathbf{\tilde{\Sigma}},Mn-p).\vspace{-.7cm}\]
			\end{proof}
			
			\begin{proof}[Proof of Corollary \ref{cor:dist}]
				
				The proof is identical to that
of Theorem \ref{thm:dist} replacing $\overline{\mathbf{S}}^{\bullet}_M$ by $\mathbf{S}^\bullet_{comb}$, noting that from Corollary 
				\ref{cor:pdf}, conditional on $\tilde{\mathbf{\Sigma}}$, $\mathbf{\overline{B}}^{\bullet}_M$ is  $N_{pm}(\mathbf{B,(\Sigma+\frac{1}{M}\tilde{\Sigma})\otimes(XX')^{-1}})$  
				and  $(Mn-p)\mathbf{S}^\bullet_{comb}$ is $W_m(\mathbf{\tilde{\Sigma}},Mn-p)$, independent of $\mathbf{\overline{B}}^{\bullet}_M$.
			\end{proof}

			\section{Details on several results}\label{App:last}

			\subsection{The posterior distributions for $\mathbf{\Sigma}$ and $\mathbf{B}$}\label{Aapp:last0}
			
Let us start by observing that $\mathbf{Y|_{B,\Sigma}}\sim N_{mn}(\mathbf{B'X,I_n\otimes \Sigma})$ and that the likelihood function for $\mathbf{Y}$ will be 
			$$l(\mathbf{B,\Sigma}|_{\mathcal{Y}})\propto |\mathbf{\Sigma}|^{-n/2} e^{-\frac{1}{2}tr\lbrace\mathbf{\Sigma}^{-1}(\mathbf{Y-B'X})(\mathbf{Y-B'X})'\rbrace}.$$
			We may then get the joint posterior distribution of $(\mathbf{B,\Sigma})$ from the product of the prior and likelihood functions as
			\begin{equation}\label{eq:posterior}
			\pi(\mathbf{B,\Sigma}|_{\mathcal{Y}})\propto |\mathbf{\Sigma}|^{-\frac{n+\alpha}{2}} e^{-\frac{1}{2}tr\lbrace\mathbf{\Sigma}^{-1}(\mathbf{Y-B'X})(\mathbf{Y-B'X})'\rbrace}.
			\end{equation}
			
			The exponent in (\ref{eq:posterior}) may be written as 
			{\small \begin{align*}
				&tr\lbrace\mathbf{\Sigma}^{-1}(\mathbf{Y-B'X})(\mathbf{Y-B'X})'\rbrace=tr\big\{\mathbf{\Sigma}^{-1}(\mathbf{Y-\hat{B}'X+\hat{B}'X-B'X})\\
				&\hspace{240pt}\times(\mathbf{Y-\hat{B}'X+\hat{B}'X-B'X})'\big\}\\
				&=tr\left\{\mathbf{\Sigma}^{-1}\left[(\mathbf{Y-\hat{B}'X})(\mathbf{Y-\hat{B}'X})'\right]\right\}\\
				&\hspace{40pt}+tr\Big\{\mathbf{\Sigma}^{-1}\big[(\mathbf{Y-\hat{B}'X})(\mathbf{\hat{B}'X-B'X})'+(\mathbf{\hat{B}'X-B'X})(\mathbf{Y-\hat{B}'X})'\\
				&\hspace{230pt}+(\mathbf{\hat{B}'X-B'X})(\mathbf{\hat{B}'X-B'X})'\big]\Big\}\\
				&=tr\left\{\mathbf{\Sigma}^{-1}\left[(\mathbf{Y-\hat{B}'X})(\mathbf{Y-\hat{B}'X})'\right]+(\mathbf{B}-\mathbf{\hat{B}})'(\mathbf{XX'})(\mathbf{B}-\mathbf{\hat{B}})\right\}\\
				& \hspace{190pt}+2 tr\left\{{\mathbf{\Sigma}}^{-1}\left[(\mathbf{Y-\hat{B}'X})(\mathbf{\hat{B}'X-B'X})'\right]\right\},
				\end{align*}}
			where, using $\mathbf{\hat{B}}'=\left[(\mathbf{XX'})^{-1}\mathbf{XY'}\right]'=\mathbf{YX'(XX')^{-1}}$,
			\begin{align*}
			(\mathbf{Y-\hat{B}'X})(\mathbf{\hat{B}'X-B'X})' & =\mathbf{YX'\hat{B}-YX'B+\hat{B}XX'\hat{B}+\hat{B}XX'B}\\
			&= \mathbf{YX'\hat{B}-YX'B+\mathbf{YX'(XX')}^{-1}}\mathbf{XX'\hat{B}}\\
			&\hspace{120pt}	+\mathbf{YX'(XX')}^{-1}\mathbf{XX'B} \\
			&=\mathbf{YX'\hat{B}-YX'B}-\mathbf{YX'\hat{B}+YX'B}=0.
			\end{align*}
			Therefore, the joint posterior distribution of $(\mathbf{B},\mathbf{\Sigma})$ is proportional to
			\begin{flalign*}
			|\mathbf{\Sigma}|^{-\frac{n+\alpha-p}{2}} e^{-\frac{n-p}{2}tr\lbrace\mathbf{\Sigma}^{-1}\mathbf{S}\rbrace}\times|\mathbf{\Sigma}|^{-\frac{p}{2}} e^{-\frac{1}{2}tr\lbrace\mathbf{\Sigma}^{-1}(\mathbf{B}-\mathbf{\hat{B}})'(\mathbf{XX'})(\mathbf{B}-\mathbf{\hat{B}})\rbrace}
			\end{flalign*}
			
			In conclusion, by Corollary 2.4.6.2. in \cite{kollo05}, the posterior distribution for $\mathbf{\Sigma}$ is 
			$$\mathbf{\Sigma}|_{\mathbf{S}}\sim W_m^{-1} \left((n\!-\!p)\mathbf{S},\!n+\!\alpha\!-p\right)\implies \mathbf{\Sigma}^{-1}|_{\mathbf{S}}\sim W_m \left(\frac{1}{n\!-\!p}\mathbf{S}^{-1},n\!+\!\alpha\!-\!p\!-\!m\!-\!1\right)$$
			and the posterior distribution for $\mathbf{B}$ is
			$$\mathbf{B}|_{\hat{\mathbf{B}},\mathbf{\Sigma}}\sim N_{pm}(\mathbf{\hat{B},\Sigma\otimes(XX')^{-1}}),$$
			assuming $n+\alpha>p+m+1$.

\vspace{.35cm}

			\subsection{Matrix calculations required in the proof of Theorem \ref{thm:pdf}}\label{Aapp:last2}
			
			
			For $\mathbf{\tilde{B}}$, $\mathbf{B}$ and $\mathbf{X}$ defined as in Section \ref{sec:post} we have
			\begin{align*} M(\mathbf{\tilde{B}}-\mathbf{\overline{B}}^\bullet_M)'\mathbf{XX'}(\mathbf{\tilde{B}}-\mathbf{\overline{B}}^\bullet_M)+(\mathbf{\tilde{B}-\hat{B}})'\mathbf{XX'}(\mathbf{\tilde{B}-\hat{B}}) = & \\
			& \hspace{-250pt}= (M+1)\mathbf{\tilde{B'}}\mathbf{XX'}\mathbf{\tilde{B}}-M\mathbf{\overline{B}}^{\bullet'}_M\mathbf{XX'}\mathbf{\tilde{B}}-M\mathbf{\tilde{B}}'\mathbf{XX'}\mathbf{\overline{B}}^\bullet_M +M\mathbf{\overline{B}}^{\bullet'}_M\mathbf{XX'}\mathbf{\overline{B}}^\bullet_M\\
			& \hspace{-60pt} -\mathbf{\hat{B}'}\mathbf{XX'}\mathbf{\tilde{B}}-\mathbf{\tilde{B}}'\mathbf{XX'}\mathbf{\hat{B}}+\mathbf{\hat{B}}'\mathbf{XX'}\mathbf{\hat{B}}\\
			&\hspace{-250pt}= (M+1)\mathbf{\tilde{B'}}\mathbf{XX'}\mathbf{\tilde{B}}-\mathbf{\tilde{B}}'\mathbf{XX'}(M\mathbf{\overline{B}}^\bullet_M+\mathbf{\hat{B})}-(M\mathbf{\overline{B}}^\bullet_M+\mathbf{\hat{B})}'\mathbf{XX'}\mathbf{\tilde{B}}\\
			& \hspace{-40pt}
			+M\mathbf{\overline{B}}^{\bullet'}_M\mathbf{XX'}\mathbf{\overline{B}}^\bullet_M+\mathbf{\hat{B}}'\mathbf{XX'}\mathbf{\hat{B}}\\
			& \hspace{-250pt}= (M+1)\left[\mathbf{\tilde{B}}-\frac{1}{M+1}(M\mathbf{\overline{B}}^\bullet_M+\mathbf{\hat{B})}\right]'\mathbf{XX'}\left[\mathbf{\tilde{B}}-\frac{1}{M+1}(M\mathbf{\overline{B}}^\bullet_M+\mathbf{\hat{B})}\right]\\
			& \hspace{-250pt} ~~~~~~~~+M\mathbf{\overline{B}}^{\bullet'}_M\mathbf{XX'}\mathbf{\overline{B}}^\bullet_M+\mathbf{\hat{B}}'\mathbf{XX'}\mathbf{\hat{B}}-\frac{1}{M+1}(M\mathbf{\overline{B}}^\bullet_M+\mathbf{\hat{B}})'\mathbf{XX}'(M\mathbf{\overline{B}}^\bullet_M+\mathbf{\hat{B}}).
			\end{align*}
			
			Since, 
			\begin{equation*}
			\begin{aligned}
			M\mathbf{\overline{B}}^{\bullet'}_M\mathbf{XX'}\mathbf{\overline{B}}^\bullet_M+\mathbf{\hat{B}}'\mathbf{XX'}\mathbf{\hat{B}}-\frac{1}{M+1}(M\mathbf{\overline{B}}^\bullet_M+\mathbf{\hat{B}})'\mathbf{XX}'(M\mathbf{\overline{B}}^\bullet_M+\mathbf{\hat{B}}) &\\
			& \hspace{-310pt}=M\mathbf{\overline{B}}^{\bullet'}_M\mathbf{XX'}\mathbf{\overline{B}}^\bullet_M+\mathbf{\hat{B}}'\mathbf{XX'}\mathbf{\hat{B}} \\
			& \hspace{-220pt} -\frac{M^2}{M+1}\mathbf{\overline{B}}^{\bullet'}_M\mathbf{XX'}\mathbf{\overline{B}}^{\bullet}_M-\frac{1}{M+1}\mathbf{\hat{B}}'\mathbf{XX'}\mathbf{\hat{B}}\\
			&\hspace{-140pt} -\frac{M}{M+1}\mathbf{\overline{B}}^{\bullet'}_M\mathbf{XX'}\mathbf{\hat{B}}-\frac{M}{M+1}\mathbf{\hat{B}}'\mathbf{XX'}\mathbf{\overline{B}}^{\bullet}_M \\
			& \hspace{-310pt} =\frac{M}{M+1}\mathbf{\overline{B}}^{\bullet'}_M\mathbf{XX'}\mathbf{\overline{B}}^{\bullet}_M+\frac{M}{M+1}\mathbf{\hat{B}}'\mathbf{XX'}\mathbf{\hat{B}}-\frac{M}{M+1}\mathbf{\overline{B}}^{\bullet'}_M\mathbf{XX'}\mathbf{\hat{B}}\\
			&\hspace{-45pt}
			-\frac{M}{M+1}\mathbf{\hat{B}}'\mathbf{XX'}\mathbf{\overline{B}}^{\bullet}_M\\
			& \hspace{-310pt}=\frac{M}{M+1}(\mathbf{\overline{B}}^{\bullet}_M-\mathbf{\hat{B}})'\mathbf{XX}'(\mathbf{\overline{B}}^{\bullet}_M-\mathbf{\hat{B}})
			\end{aligned}
			\end{equation*}
			we may write
			\begin{align*}
			&M(\mathbf{\tilde{B}}-\mathbf{\overline{B}}^\bullet_M)'\mathbf{XX'}(\mathbf{\tilde{B}}-\mathbf{\overline{B}}^\bullet_M)+(\mathbf{\tilde{B}-\hat{B}})'\mathbf{XX'}(\mathbf{\tilde{B}-\hat{B}})=\\ &=(M+1)\left[\mathbf{\tilde{B}}-\frac{1}{M+1}(M\mathbf{\overline{B}}^\bullet_M+\mathbf{\hat{B}})\right]'\mathbf{XX'}\left[\mathbf{\tilde{B}}-\frac{1}{M+1}(M\mathbf{\overline{B}}^\bullet_M+\mathbf{\hat{B}})\right] &\\
			& \hspace{200pt} +\frac{M}{M+1}(\mathbf{\overline{B}}^\bullet_M-\mathbf{\hat{B}})'\mathbf{XX}'(\mathbf{\overline{B}}^\bullet_M-\mathbf{\hat{B}})
			\end{align*}

			\subsection{Details about the derivations of results 1, 2 and 5 in Section \ref{ssec:mul1}}\label{Aapp:last3}
			
			
			\noindent\textit{Details on Result 1} 
			
			From (\ref{eq:likelihood}) we may immediately conclude that the MLE of $\mathbf{B}$ based on the synthetic data will be $\mathbf{\overline{B}}^{\bullet}_M$ with
			$$E(\mathbf{\overline{B}}^{\bullet}_M)=(\mathbf{XX}')^{-1}\mathbf{X}\frac{1}{M}\sum_{j=1}^{M}E(\mathbf{W}'_j)=(\mathbf{XX}')^{-1}\mathbf{XX'}E(\mathbf{\tilde{B}})=E(\mathbf{\hat{B}})=\mathbf{B}$$
			and
			\begin{equation}\label{eq:var1}
			Var(\mathbf{\overline{B}}^{\bullet}_M)=Var[E(\mathbf{\overline{B}}^{\bullet}_M|\mathbf{\tilde{B},\tilde{\Sigma}})]+E[Var(\mathbf{\overline{B}}^{\bullet}_M|\mathbf{\tilde{B},\tilde{\Sigma}})].
			\end{equation}
			
			For the first term in (\ref{eq:var1}), we have
			$$Var[E(\mathbf{\overline{B}}^{\bullet}_M|\mathbf{\tilde{B},\tilde{\Sigma}})]=Var[\mathbf{\tilde{B}}]=Var[E(\mathbf{\tilde{B}|\hat{B},\tilde{\Sigma}})]+E[Var(\mathbf{\tilde{B}|\hat{B},\tilde{\Sigma}})]=$$
			$$=Var(\mathbf{\hat{B}})+E[\mathbf{\tilde{\Sigma}\otimes (XX')^{-1}}]=\mathbf{\Sigma\otimes (XX')^{-1}}+\frac{n-p}{n+\alpha-p-2m-2}\mathbf{\Sigma\otimes (XX')^{-1}}$$
			and for the second term, we have
			$$E[Var(\mathbf{\overline{B}}^{\bullet}_M|\mathbf{\tilde{B},\tilde{\Sigma}})]=E\left[\frac{1}{M}\mathbf{\tilde{\Sigma}\otimes (XX')^{-1}}\right]=\frac{1}{M}\frac{n-p}{n+\alpha-p-2m-2}\mathbf{\Sigma\otimes (XX')^{-1}},$$
			so that 
			$$Var(\mathbf{\overline{B}}^{\bullet}_M)=\frac{2M(n-p-m-1)+n-p+M\alpha}{M(n+\alpha-p-2m-2)}\mathbf{\Sigma\otimes (XX')^{-1}}$$
			under the condition that $n+\alpha>p+2m+2$.

			\noindent \textit{Details on Result 2} $$E(\mathbf{\overline{S}}^{\bullet}_M)=E(\mathbf{\tilde{\Sigma}})=E\left(\frac{n-p}{n+\alpha-p-2m-2}\mathbf{S}\right)=\frac{n-p}{n+\alpha-p-2m-2}\mathbf{\Sigma}.$$
			
\pagebreak
			
			\noindent \textit{Details on Result 5} 
			
			Let us consider $\mathbf{H}$ and $\mathbf{G}$ given by (\ref{eq:H}) and (\ref{eq:G}). We will begin by rewriting all four classical statistics $T_{1,M}^\bullet$, $T_{2,M}^\bullet$, $T_{3,M}^\bullet$ and $T_{4,M}^\bullet$ in Subsection \ref{ssec:mul1}, in order to make them assume the same kind of form and then we will prove why all of them are non-pivotal, without loss of generality considering $M=1$. The first statistic, $T_{1,M}^\bullet\,$may be rewritten as
			\begin{flalign*}
			T^\bullet_{1,1}   & =  \frac{|\mathbf{G}|}{|\mathbf{G}+(n-p)\mathbf{\tilde{\Sigma}^{-1/2}(2\tilde{\Sigma}+\Sigma)^{1/2}H(2\tilde{\Sigma}+\Sigma)^{1/2}\tilde{\Sigma}^{-1/2}}|}.
			\end{flalign*}
			while $T_{2,M}^\bullet$ and $T_{3,M}^\bullet$ may be rewritten as
			\begin{flalign*}
			T^\bullet_{2,1} & =(n-p)tr\left[\mathbf{H(2\tilde{\Sigma}+\Sigma)^{1/2}\tilde{\Sigma}^{-1/2}}\mathbf{G}^{-1}\mathbf{\tilde{\Sigma}^{-1/2}(2\tilde{\Sigma}+\Sigma)^{1/2}}\right],
			\end{flalign*}
			\begin{flalign*}
			T^\bullet_{3,1} &=tr\lbrace\mathbf{ H \times [H+(2\tilde{\Sigma}+\Sigma)^{-1/2}\tilde{\Sigma}^{1/2}}\times(n-p)\mathbf{G\times \tilde{\Sigma}^{1/2}(2\tilde{\Sigma}+\Sigma)^{-1/2}]^{-1}}\rbrace.
			\end{flalign*}
			Concerning $T^\bullet_{4,1}$, we have $T^\bullet_{4,1}=\lambda_1$ where $\lambda_1$ denotes the largest eigenvalue of 
			\begin{flalign*}
			&(n-p)\mathbf{H\times(2\tilde{\Sigma}+\Sigma)^{1/2}\tilde{\Sigma}^{-1/2}\times G^{-1}\times \tilde{\Sigma}^{-1/2}(2\tilde{\Sigma}+\Sigma)^{1/2}}.
			\end{flalign*}
			
			We can observe that a term in the denominator of the expression $T^\bullet_{1,1}$ is
			$$\mathbf{\tilde{\Sigma}^{-1/2}(2\tilde{\Sigma}+\Sigma)^{1/2}H(2\tilde{\Sigma}+\Sigma)^{1/2}\tilde{\Sigma}^{-1/2}}|_\mathbf{\tilde{\Sigma}^{-1}}\sim W_m(\mathbf{(2I+\tilde{\Sigma}^{-1/2}\Sigma \tilde{\Sigma}^{-1/2})},p),$$ 
			while in the expressions for the other statistics there are similar terms. These terms involve a product similar to $\mathbf{\tilde{\Sigma}}^{-1/2}(\mathbf{2\tilde{\Sigma}+\Sigma})^{1/2}$ that cannot be simplified to an expression which is not a function of $\mathbf{\Sigma}$, therefore making these statistics non-pivotal.
			\noindent

\begin{figure}[h!]
				
				\begin{subfigure}{.5\textwidth}
					\centering
					\includegraphics[width=.7\linewidth]{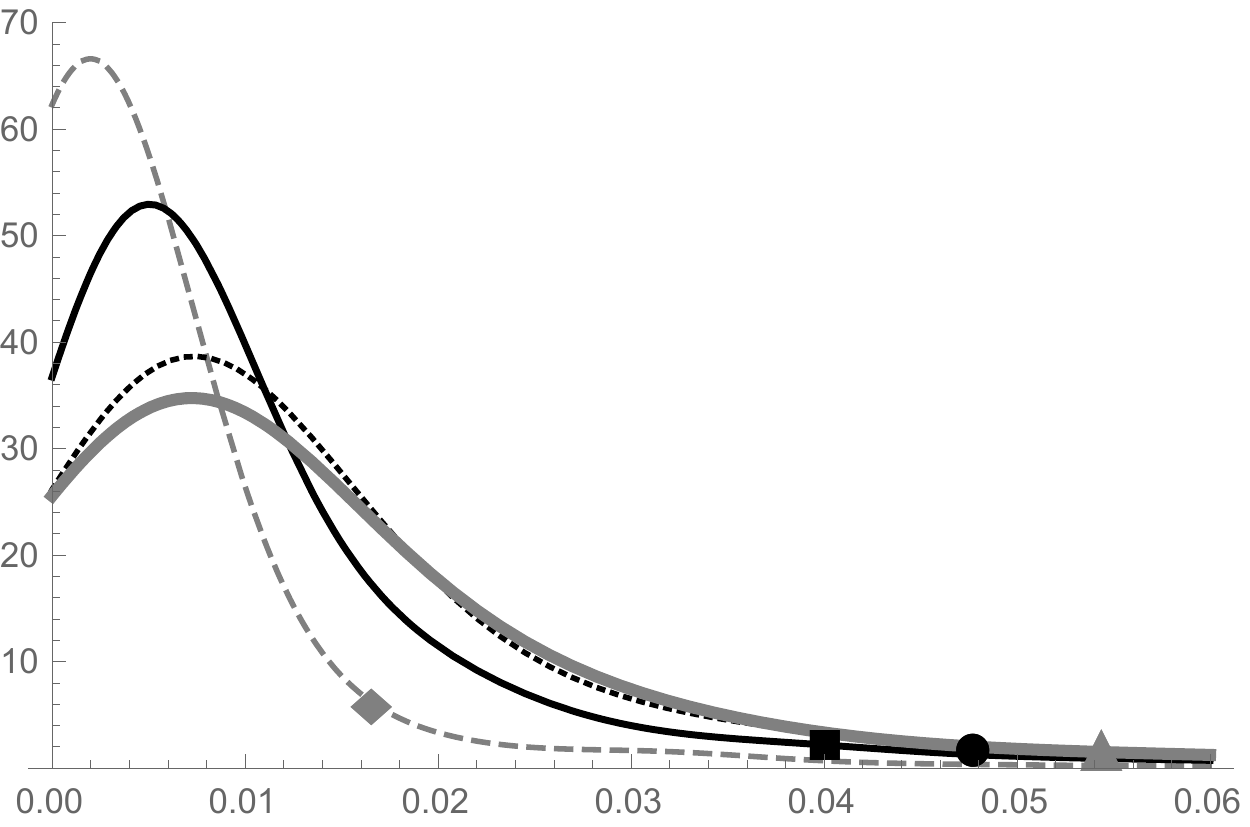}\includegraphics[width=.15\linewidth]{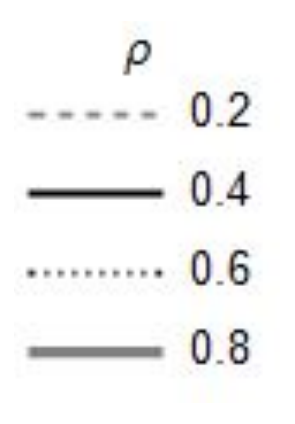}
					\label{fig:wilks}
					\vspace{-5pt}
					\caption{Wilks}
				\end{subfigure}%
				\vspace{-6pt}
				\begin{subfigure}{.5\textwidth}
					\centering
					\includegraphics[width=.7\linewidth]{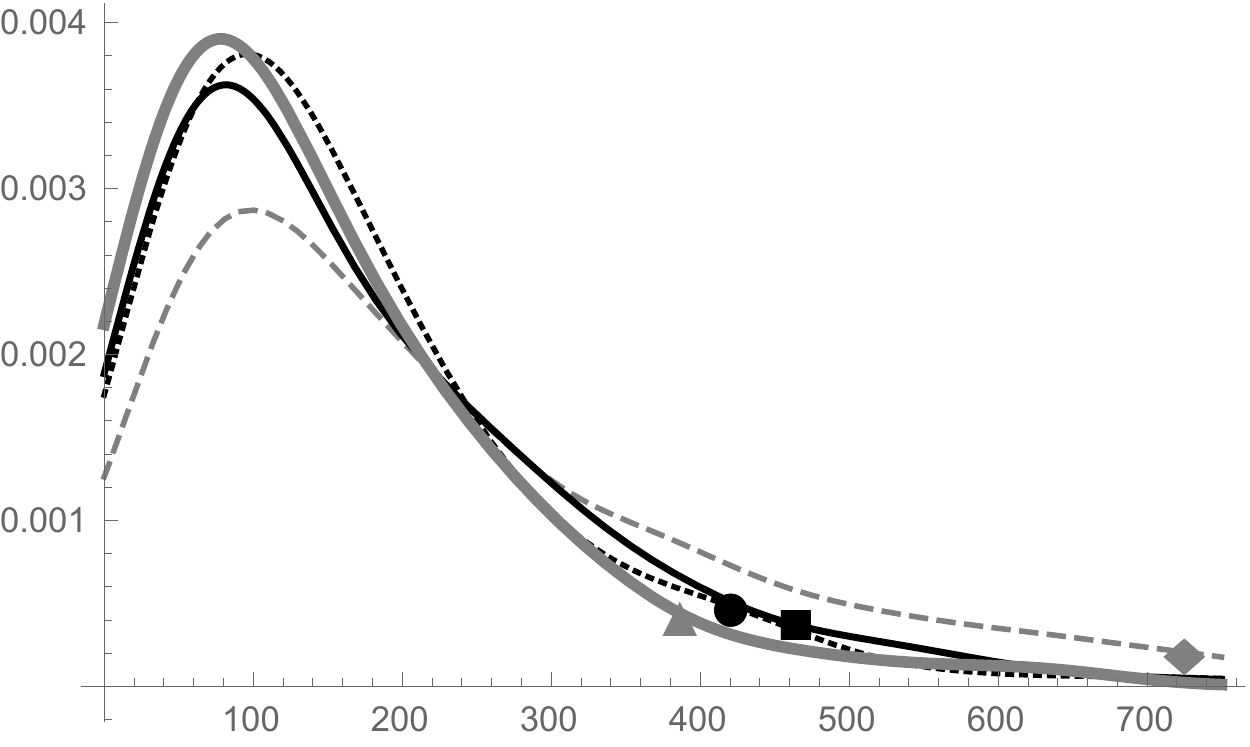}\includegraphics[width=.15\linewidth]{legend_bw_rho}
					\label{fig:lawley}
					\vspace{-5pt}
					\caption{Lawley}
				\end{subfigure}%
				\vspace{-6pt}
				
				\begin{subfigure}{.5\textwidth}
					\centering
					\includegraphics[width=.7\linewidth]{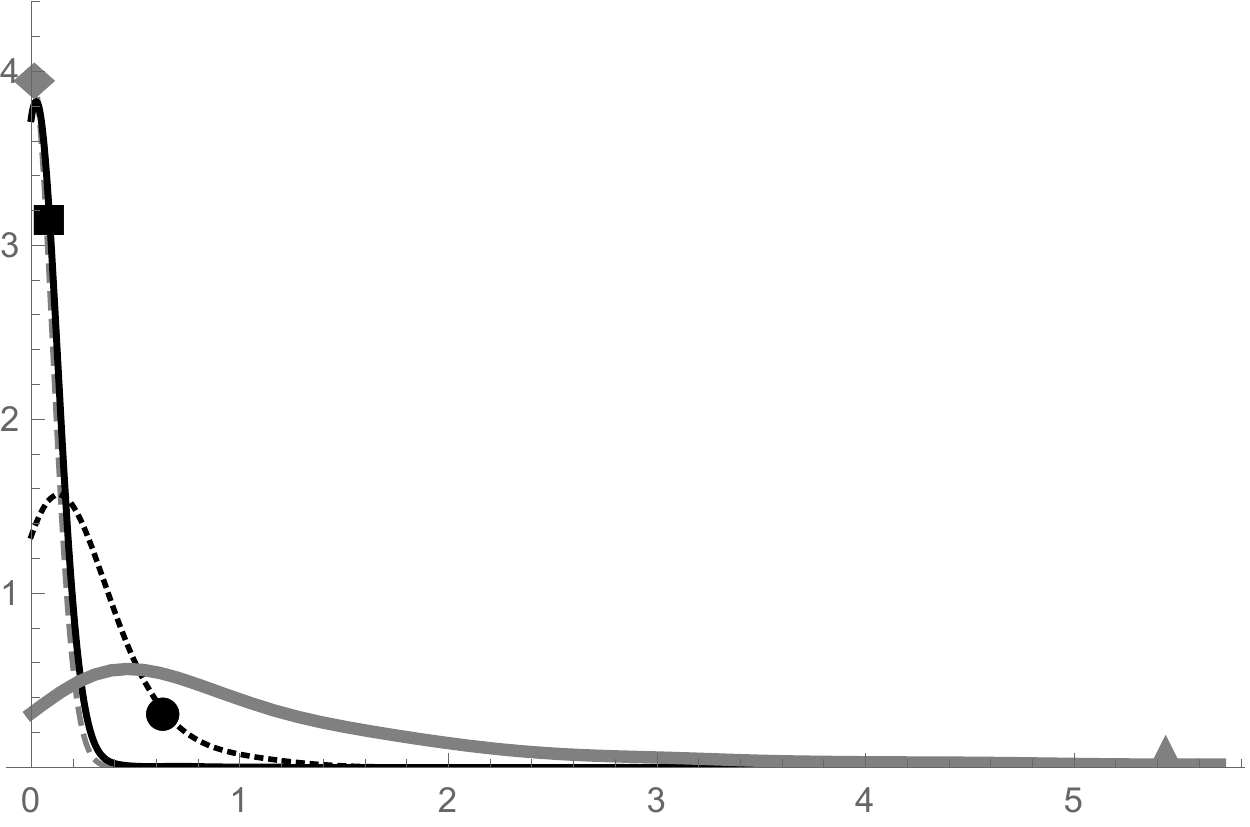}\includegraphics[width=.15\linewidth]{legend_bw_rho}
					\label{fig:pillai}
					\vspace{-5pt}
					\caption{Pillai}
				\end{subfigure}%
				\vspace{-6pt}
				\begin{subfigure}{.5\textwidth}
					\centering
					\includegraphics[width=.7\linewidth]{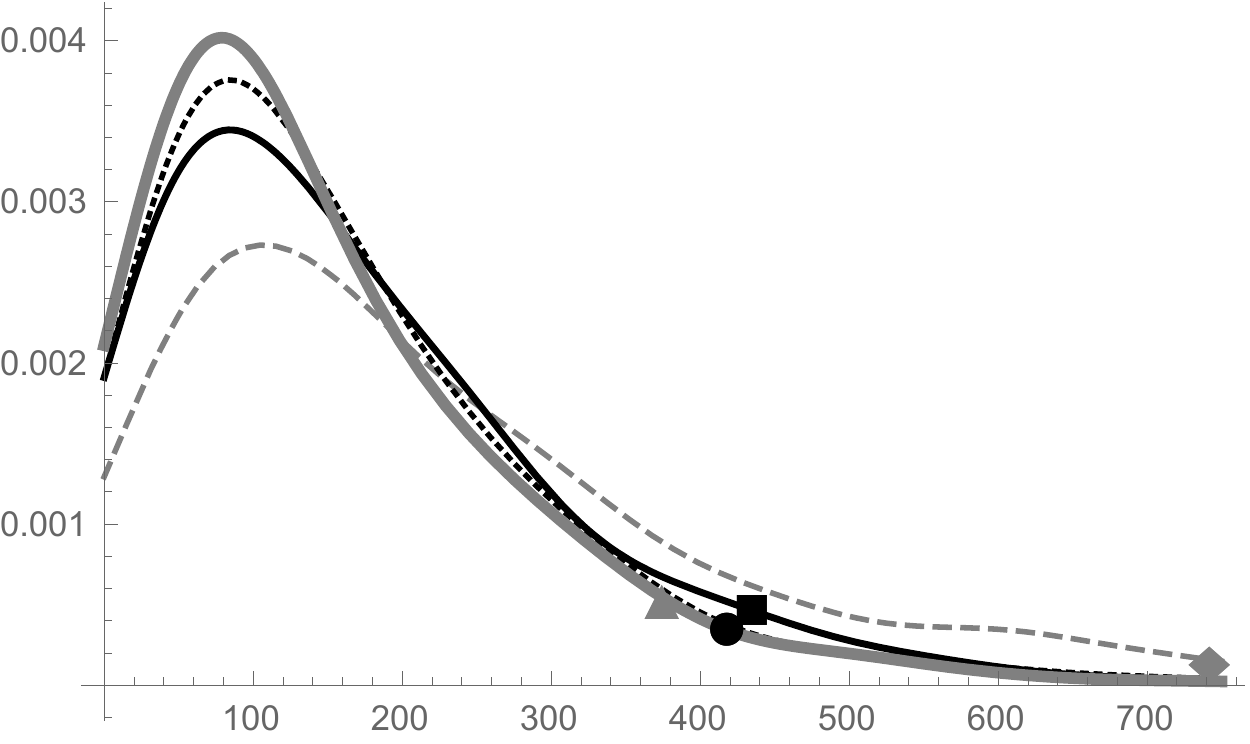}\includegraphics[width=.15\linewidth]{legend_bw_rho}
					\label{fig:roy}
					\vspace{-5pt}
					\caption{Roy}
				\end{subfigure}%
				\vspace{-5pt}
				\caption{Smoothed empirical\! distributions\! and\! cut-off\! points (${\gamma\!=\!0.05}$) of $T^\bullet_{1,1}$, $T^\bullet_{2,1}$, $T^\bullet_{3,1}$ and $T^\bullet_{4,1}$ for $\rho=\lbrace\!0.2,\! 0.4,\! 0.6,\! 0.8\rbrace$.}\label{fig:fig}
				
			\end{figure}

			Thus, in order to illustrate how these statistics are dependent on $\mathbf{\Sigma}$, we can analyze in Figure \ref{fig:fig} the empirical distributions of $T^\bullet_{1,1}$, $T^\bullet_{2,1}$, $T^\bullet_{3,1}$ and $T^\bullet_{4,1}$ when we consider a simple case where $m=2$, $p=3$, $\alpha=4$, $n=100$ and $\Sigma=\left( \begin{smallmatrix} 1&\rho\\ \rho&1 \end{smallmatrix} \right)$ with $\rho=\lbrace0.2, 0.4, 0.6, 0.8\rbrace$ for a simulation size of 1000.

			\subsection{Details about the derivation of result 1 in Subsection \ref{ssec:mul2}}\label{Aapp:last5}
			
			
			
			Recalling that $(Mn-p) \mathbf{S}^\bullet_{comb}|_\mathbf{\tilde{\Sigma}}\sim W_m(\mathbf{\tilde{\Sigma}},Mn-p)$ and 
			that $\mathbf{\tilde{\Sigma}}^{-1}|_{\mathbf{S}}\sim W_m (\frac{1}{n-p}\mathbf{S}^{-1},n+\alpha-p-m-1)$ we immediately obtain 
			$$E(\mathbf{S}^{\bullet}_{comb})=E(\mathbf{\tilde{\Sigma}})=E\left(\frac{n-p}{n+\alpha-p-2m-2}\mathbf{S}\right)=\frac{n-p}{n+\alpha-p-2m-2}\mathbf{\Sigma}.$$

			\subsection{Details about the derivations of the results in Section \ref{sec:sim}}\label{Aapp:last6}

			\noindent \textit{Details on the Expected Values in Section \ref{sec:sim}} 
			
			Recall that $(n-p)\mathbf{S}\sim W_m(\mathbf{\Sigma},n-p)$, 
			thus implying that $$E(|(n-p)\mathbf{S}|)=|\mathbf{\Sigma}|E(\prod_{i=1}^m\chi^2_{n-p-i+1})=\frac{(n-p)!}{(n-p-m)!}|\mathbf{\Sigma}|,$$ and recall that 
			$$\mathbf{\tilde{\Sigma}}|_\mathbf{S}\sim W^{-1}_m((n-p)\mathbf{S},n+\alpha-p)\implies \mathbf{\tilde{\Sigma}}^{-1}|_\mathbf{S}\sim W_m\left(\frac{1}{n-p}\mathbf{S}^{-1},n\!+\!\alpha\!-\!p\!-\!m\!-\!1\right)$$
			thus implying that, making $\kappa_{n,\alpha,p,m}=n+\alpha-p-m-1$, given $\mathbf{S}$,
			\begin{flalign*}
			E(|\mathbf{\tilde{\Sigma}}|) & =E(|\mathbf{\tilde{\Sigma}}^{-1}|^{-1})=|(n-p)\mathbf{S}|E\left(\frac{1}{\prod_{i=1}^m \chi^2_{\kappa_{n,\alpha,p,m}-i+1}}\right)\\
			&=|(n-p)\mathbf{S}|\frac{(-2+\kappa_{n,\alpha,p,m}-m)!}{(-2+\kappa_{n,\alpha,p,m})!},
			\end{flalign*}
			since $\prod_{i=1}^m \chi^2_{\kappa_{n,\alpha,p,m}-i+1}$ is a product of independent $\chi^2$ variables. Also recalling that, given $\mathbf{\tilde{\Sigma}}$,  we have $M(n-p)\mathbf{\overline{S}^{\bullet}}_M\sim W_m(\mathbf{\tilde{\Sigma}},M(n-p))$ 
			and $(Mn-p)\mathbf{S}^\bullet_{comb}\sim W_m(\mathbf{\tilde{\Sigma}},Mn-p)$, we may conclude that, given $\mathbf{\tilde{\Sigma}}$, 
			$$E(|M(n-p)\mathbf{S}^{\bullet}_M|)=\frac{(Mn-Mp)!}{(Mn-Mp-m)!}\times|\mathbf{\tilde{\Sigma}}|$$ and 
			$$E(|(Mn-p)\mathbf{S}^\bullet_{comb}|)=\frac{(Mn-p)!}{(Mn-p-m)!}\times|\mathbf{\tilde{\Sigma}}|.$$

			Combining the results for $E(|(n-p)\mathbf{S}|)$ and $E(|\mathbf{\tilde{\Sigma}}|)|_\mathbf{S}$  with each of the expected values for $|M(n-p)\mathbf{S}^{\bullet}_M|$ and $|(Mn-p)\mathbf{S}^\bullet_{comb}|$, we end up with the expression 
			for $E(\!\Upsilon_M\!)$ found in Section \ref{sec:sim}.

				\section{Joining multiple datasets into a single dataset}\label{App:C}
				
				Let us consider the $M$ synthetic datasets as one only dataset of size $nM$
				$$
				\left(\ba{c}
				{\mathbf W}_a\\
				{\mathbf X}_a
				\ea\right)=
				\left(~~\begin{array}{c|c|c|c}
				\mathbf{W}_1  & \mathbf{W}_2 & \dots & \mathbf{W}_M \\
				\hline
				\mathbf{X} &\mathbf{X} & \dots &\mathbf{X}
				\end{array}~~\right),
				$$
				where $\mathbf{W}_a=(\mathbf{W}_1| ... | \mathbf{W}_M)$ is the $m\times nM$ matrix of the synthesized data under FPPS and $\mathbf{X}_a=(\mathbf{X}|...|\mathbf{X})$ the $p\times nM$ matrix of the $M$ repeated `fixed' sets of covariates, from the original data.

				Let
				$$
				\mathbf{B}_a=(\mathbf{X}_a \mathbf{X}_a')^{-1}\mathbf{X}_a \mathbf{W}_a'
				$$
				be the estimator for $\mathbf{B}$, based on the dataset of size $nM$, obtained by joining the $M$ synthetic datasets in one only dataset. Consequently one has that
				$$
				\begin{array}{rcl}
				\mathbf{B}_a & = & \ds (\mathbf{X}_a \mathbf{X}_a')^{-1}\mathbf{X}_a \mathbf{W}_a'=(M (\mathbf{X} \mathbf{X}'))^{-1}\mathbf{X}_a \mathbf{W}_a'=\frac{1}{M}( \mathbf{X} \mathbf{X}')^{-1}\mathbf{X}_a \mathbf{W}_a'\ms\\
				& = & \ds \frac{1}{M}(\mathbf{X} \mathbf{X}')^{-1} \Bigl(\underbrace{\mathbf{X}| ~\dots~ |\mathbf{X}}_{M\; times}\Bigr) \mathbf{W}_a'=\frac{1}{M}\left(( \mathbf{X} \mathbf{X}')^{-1} \mathbf{X} \mathbf{W}_1 \!+\!\dots\!+\! (  \mathbf{X} \mathbf{X}')^{-1} \mathbf{X}\mathbf{W}_M\right)\ms\\
				& = & \ds \frac{1}{M}(\mathbf{X} \mathbf{X}')^{-1} \mathbf{X} \left( \mathbf{W}_1 +...+\mathbf{W}_M\right)=(\mathbf{X} \mathbf{X}')^{-1} \mathbf{X}\mathbf{\overline{W}}_{M} =\mathbf{\overline{B}}_M^\bullet\,,
				\ea
				$$
				which is same estimator for $\mathbf{B}$ as in (\ref{eq:par2nd1}).
				
				Now let
				$$
				\mathbf{S}_a=\frac{1}{nM-p}(\mathbf{W}_a-\mathbf{B}_a'\mathbf{X}_a)(\mathbf{W}_a-\mathbf{B}_a'\mathbf{X}_a)'
				$$
				be the estimator for $\mathbf{\Sigma}$, based on the dataset of size $nM$, obtained by joining the $M$ synthetic datasets in one only dataset.
				
				Observe that $\mathbf{\overline{W}}_{M}=\frac{1}{M}\sum_{j=1}^{M}\mathbf{W}_j$, defined before expression (\ref{eq:par2nd1}), can be written as
				$$
				\mathbf{\overline{W}}_{M}=\frac{1}{M}\mathbf{W}_a \mathbf{R}
				$$
				with $\mathbf{R}=\left(\overrightarrow{\mathbf{1}}_M\otimes \mathbf{I}_n\right)$ where $\overrightarrow{\mathbf{1}}_M$ is a vector of $1$'s of size $M$. 
				
				Now let us consider the estimator $\mathbf{S}_\mathbf{w}$ of $\mathbf{\Sigma}$, defined in the text, before expression (\ref{eq:par2nd1}). This estimator may be written as
				$$
				\mathbf{S}_\mathbf{w}=\sum_{i=1}^n\sum_{j=1}^M (\mathbf{w}_{ji}-\mathbf{\overline{w}}_i)(\mathbf{w}_{ji}-\mathbf{\overline{w}}_i)',
				$$
				where $\mathbf{w}_{ji}$ is the $i$-th column of $\mathbf W_j$ $(i=1,\dots,n;j=1,\dots,M)$. We may thus write
				$$
				\ba{rcl}
				\mathbf{S}_\mathbf{w} & = & \ds \left(\mathbf{W}_a-\overrightarrow{\mathbf{1}}_M'\otimes\mathbf{\overline{W}}_{M}\right)\left(\mathbf{W}_a-\overrightarrow{\mathbf{1}}_M'\otimes\mathbf{\overline{W}}_{M}\right)'\ms\\
				& = & \ds \left(\mathbf{W}_a-\frac{1}{M}\overrightarrow{\mathbf{1}}_M'\otimes(\mathbf{W}_a \mathbf{R})\right)\left(\mathbf{W}_a-\frac{1}{M}\overrightarrow{\mathbf{1}}_M'\otimes(\mathbf{W}_a \mathbf{R})\right)'\ms\\
				& = & \ds \left(\mathbf{W}_a-\frac{1}{M}\mathbf{W}_a\mathbf{R}\mathbf{R}'\right)\left(\mathbf{W}_a-\frac{1}{M}\mathbf{W}_a\mathbf{R}\mathbf{R}'\right)'
				\ea
				$$
				and the estimator $\mathbf{S}_{mean}$ of $\mathbf{\Sigma}$, defined right after expression (\ref{eq:par2nd2}) as
				$$
				\mathbf{S}_{mean}=\left(\frac{1}{M}\mathbf{W}_a\mathbf{R}-\frac{1}{M}\mathbf{B}_a'\mathbf{X}_a\mathbf{R}\right)\left(\frac{1}{M}\mathbf{W}_a\mathbf{R}-\frac{1}{M}\mathbf{B}_a'\mathbf{X}_a\mathbf{R}\right)'.$$
				
				We may therefore write the combination estimator $\mathbf{S}_{comb}$ defined in (\ref{eq:par2nd2}) as
				$$
				\ba{l} \ds \mathbf{S}_{comb}=\frac{1}{nM-p}\left[\left(\mathbf{W}_a-\frac{1}{M}\mathbf{W}_a\mathbf{R}\mathbf{R}'\right)\left(\mathbf{W}_a-\frac{1}{M}\mathbf{W}_a\mathbf{R}\mathbf{R}'\right)'\right]\ms\\
				\hskip 1.5cm \ds+\frac{1}{nM-p}\left[M\times \left(\frac{1}{M}\mathbf{W}_a\mathbf{R}-\frac{1}{M}\mathbf{B}_a'\mathbf{X}_a\mathbf{R}\right)\left(\frac{1}{M}\mathbf{W}_a\mathbf{R}-\frac{1}{M}\mathbf{B}_a'\mathbf{X}_a\mathbf{R}\right)'\right]
				\ea
				$$

				To prove that $\mathbf{S}_{comb}$ is equal to $\mathbf{S}_a$ it will only be necessary to focus on 
				$$
				\ba{l}\ds
				\left(\mathbf{W}_a-\frac{1}{M}\mathbf{W}_a\mathbf{R}\mathbf{R}'\right)\left(\mathbf{W}_a-\frac{1}{M}\mathbf{W}_a\mathbf{R}\mathbf{R}'\right)'\ms\\
				\hskip3.1cm \ds +M\times \left(\frac{1}{M}\mathbf{W}_a\mathbf{R}-\frac{1}{M}\mathbf{B}_a'\mathbf{X}_a\mathbf{R}\right)\left(\frac{1}{M}\mathbf{W}_a\mathbf{R}-\frac{1}{M}\mathbf{B}_a'\mathbf{X}_a\mathbf{R}\right)'\ms\\
				\hskip .5cm \ds =\mathbf{W}_a\mathbf{W}_a'-\frac{1}{M}\mathbf{W}_a\mathbf{R}\mathbf{R}'\mathbf{W}_a'-\frac{1}{M}\mathbf{W}_a\mathbf{R}\mathbf{R}'\mathbf{W}_a' +\frac{1}{M^2}\mathbf{W}_a\mathbf{R}\mathbf{R}'\mathbf{R}\mathbf{R}'\mathbf{W}_a'\ms\\
				\hskip6.1cm \ds +\frac{1}{M}\mathbf{W}_a\mathbf{R}\mathbf{R}'\mathbf{W}_a'-\frac{1}{M}\mathbf{B}_a'\mathbf{X}_a\mathbf{R}\mathbf{R}'\mathbf{W}_a'\\
				\hskip6.3cm \ds -\frac{1}{M}\mathbf{W}_a\mathbf{R}\mathbf{R}'\mathbf{X}_a'\mathbf{B}_a+\frac{1}{M}\mathbf{B}_a'\mathbf{X}_a\mathbf{R}\mathbf{R}'\mathbf{X}_a'\mathbf{B}_a\,,
				\ea
				$$
				which, using the fact that $\frac{1}{M}\mathbf{X}_a\mathbf{R}\mathbf{R}'=\mathbf{X}_a$ and $\frac{1}{M}\mathbf{R}\mathbf{R}'\mathbf{R}\mathbf{R}'=\mathbf{R}\mathbf{R}'$, may be written as
				$$
				\ba{l}
				\ds \mathbf{W}_a\mathbf{W}_a'-\frac{1}{M}\mathbf{W}_a\mathbf{R}\mathbf{R}'\mathbf{W}_a'-\frac{1}{M}\mathbf{W}_a\mathbf{R}\mathbf{R}'\mathbf{W}_a' +\frac{1}{M}\mathbf{W}_a\mathbf{R}\mathbf{R}'\mathbf{W}_a'\ms\\
				\hskip 3cm \ds +\frac{1}{M}\mathbf{W}_a\mathbf{R}\mathbf{R}'\mathbf{W}_a'-\mathbf{B}_a'\mathbf{X}_a\mathbf{W}_a'-\mathbf{W}_a\mathbf{X}_a'\mathbf{B}_a+\mathbf{B}_a'\mathbf{X}_a\mathbf{X}_a'\mathbf{B}_a\ms\\
				\hskip 1.5cm \ds =\mathbf{W}_a\mathbf{W}_a'-\mathbf{B}_a'\mathbf{X}_a\mathbf{W}_a'-\mathbf{W}_a\mathbf{X}_a'\mathbf{B}_a+\mathbf{B}_a'\mathbf{X}_a\mathbf{X}_a'\mathbf{B}_a\ms\\
				\hskip 1.5cm \ds =(\mathbf{W}_a-\mathbf{B}_a'\mathbf{X}_a)(\mathbf{W}_a-\mathbf{B}_a'\mathbf{X}_a)'=(nM-p)\mathbf{S}_a\,.
				\ea
				$$

		\end{appendices}
		

	\end{document}